\documentclass{amsart}
\usepackage{amsmath}
 \usepackage[colorlinks=true]{hyperref}
\hypersetup{urlcolor=blue, citecolor=blue}

\numberwithin{equation}{section} 

\newtheorem{theorem}{Theorem}[section]

\newtheorem{lemma}[theorem]{Lemma}

\theoremstyle{definition}

\newtheorem{remark}[theorem]{Remark}

\title[Almost global existence for quasi-linear wave equations]
      {On almost global existence and local well-posedness for
      some 3-D quasi-linear wave equations}

\author{Kunio Hidano}
\address{Department of Mathematics\\Faculty of Education\\Mie University\\1577 Kurima-machiya-cho, Tsu\\Mie 514-8507, Japan}
\email{hidano@edu.mie-u.ac.jp}
\thanks{The first author was partly supported by the Grant-in-Aid
for Scientific Research (C) (No.20540165), Japan Society
for the Promotion of Science.}

\author{Chengbo Wang}
\address{Department of Mathematics\\
            Johns Hopkins University\\
            Baltimore, MD 21218, USA}
\curraddr{Department of Mathematics\\
                Zhejiang University\\
                Hangzhou 310027, China}
\email{wangcbo@gmail.com}
\urladdr{http://www.math.zju.edu.cn/wang/}
\thanks{The second author was supported by the Fundamental Research Funds for the Central Universities, NSFC 10871175 and 10911120383.}

\author{Kazuyoshi Yokoyama}
\address{Hokkaido Institute of Technology, 7-15-4-1 Maeda, Teine-ku, Sapporo, Hokkaido 006-8585, Japan}\email{yokoyama@hit.ac.jp}

\dedicatory{}
\date{}

\begin{document}

\begin{abstract}
We study the Cauchy problem for a quasilinear wave equation with
 low-regularity data. A space-time
 $L^2$ estimate for the variable coefficient wave equation plays a
 central role for this purpose.
Assuming radial symmetry, we establish
 the almost global existence of a strong solution for every small initial data in $H^2
 \times H^1$. We also show that the initial value problem is locally well-posed.
\end{abstract}

\maketitle\tableofcontents


\section{Introduction}

In this paper, we consider the quasilinear wave equation
\begin{equation}
 \partial_t^2 \phi - \Delta \phi + h(\phi)\Delta \phi = F(\partial \phi)
\quad \mbox{in}\,\, (0,T)\times {\mathbb R}^3
\label{ea1}
\end{equation}
for a real-valued function $\phi = \phi(t,x)$,
obeying the initial conditions
\begin{equation}
 \phi(0,\cdot) = f,\quad
\partial_t \phi(0,\cdot) = g.
\label{ea2}
\end{equation}
We assume that the function $h$ is smooth and $h(0) = 0$.
We also assume that $F$ is smooth and quadratic with respect to
$\partial \phi = (\partial \phi/\partial t, \nabla \phi)$,
where $\nabla \phi = (\partial \phi/\partial x^1,
\partial\phi /\partial x^2, \partial \phi/\partial x^3)$.
Assuming radial symmetry, we aim at
showing the almost global existence of a low regularity solution to the initial
value problem (\ref{ea1}), (\ref{ea2}).

When the equation (\ref{ea1}) is semilinear,
namely $h \equiv 0$, Ponce-Sideris \cite{PS} proved that the
initial value problem (\ref{ea1}), (\ref{ea2}) is locally
well-posed in $H^s \times H^{s-1}$ if $s > 2$.
Smith-Tataru \cite{SmT} further proved that
the local well-posedness in $H^s \times H^{s-1}$
holds for the quasilinear case with $s > 2$.
These results are sharp in general, because
the examples given by Lindblad \cite{L1}, \cite{L2} show that
we cannot always expect the well-posedness for
quadratically nonlinear wave equations
in $H^2 \times H^1$.
However, as far as radially symmetric solutions are concerned,
there is still room for investigation.
Indeed, we know by Klainerman-Machedon \cite{KM} that
the Cauchy problem (\ref{ea1}), (\ref{ea2}) is locally
well-posed in $H^2_{{\rm rad}}\times H^1_{{\rm rad}}$, 
if the equation (\ref{ea1}) is semilinear
(here $H^m_{{\rm rad}}$ denotes the set of all functions in $H^m$ with radial
symmetry).
Furthermore, the $H^2\times H^1$-radial solutions
to the semilinear Cauchy problem
exist almost globally for small initial data
(Hidano-Yokoyama \cite{HY}).
In view of these results,
it is naturally expected to extend the study of the $H^2\times H^1$-radial solutions to the quasilinear case.
We will show that the quasilinear Cauchy problem (\ref{ea1}),
(\ref{ea2}) is locally well-posed in $H^2_{{\rm rad}} \times H^1_{{\rm rad}}$,
and the radial solutions exist almost globally for small data.

Before stating the main result,
let us introduce some function spaces.
For brevity's sake, we set $S_T = (0,T) \times {\mathbb R}^3$ in what follows.
Firstly, we set
\begin{equation}
 X_m(T) = \bigcap_{j=0}^m C^j([0,T]; H^{m-j}_{{\rm rad}}({\mathbb R}^3))
\quad \mbox{for}\,\, m = 1,2.
\nonumber
\end{equation}
For each $\phi \in X_m(T)$ we let
\[
 \|\phi\|_{X_m(T)} =
\|\phi\|_{L^{\infty}([0,T]; L^2({\mathbb R}^3))} + \|\phi\|_{E_m(T)},
\]
where
\begin{eqnarray}
 \|\phi\|_{E_1(T)} &=& \|\nabla  \phi\|_{L^{\infty}([0,T]; L^2({\mathbb R}^3))}
 + \|\partial_t \phi\|_{L^{\infty}([0,T]; L^2({\mathbb R}^3))},
 \nonumber \\
 \|\phi\|_{E_2(T)} &=& \|\phi\|_{E_1(T)}
 + \|\nabla  \phi\|_{E_1(T)}
 + \|\partial_t \phi\|_{E_1(T)}. \nonumber
\end{eqnarray}
 Secondly, we define
\begin{equation}
Y_m(T) = \{\phi\,;\, \phi \in L^1_{\rm loc}(S_T),\enskip
 \|\phi\|_{Y_m(T)} < \infty\}
\nonumber
\end{equation}
for $m = 1,2$, where
\begin{eqnarray}
 \|\phi\|^2_{Y_1(T)}\!\! &=&\!\! (1+T)^{-1/2}\left(
 \|r^{-5/4}\phi\|^2_{L^2(S_T)} +
 \|r^{-1/4}\partial \phi\|^2_{L^2(S_T)}\right),
 \nonumber \\
 \|\phi\|^2_{Y_2(T)}\!\! &=& \|\phi\|^2_{Y_1(T)}
 + \|\partial \phi\|^2_{Y_1(T)}.
 \nonumber
\end{eqnarray}
Lastly, we set
\begin{equation}
 Z_m(T) = \{\phi\,;\, \phi \in L^1_{\rm loc}(S_T),\enskip
 \|\phi\|_{Z_m(T)} < \infty\}
\nonumber
\end{equation}
for $m = 1,2$, where
\begin{eqnarray}
 \|\phi\|^2_{Z_1(T)}\!\! &=&\!\! \bigl(\log (2 + T)\bigr)^{-1}
 \bigl(\|r^{-5/4} \langle r \rangle^{-1/4} \phi\|^2_{L^2(S_T)}
\nonumber \\
 & &+
 \|r^{-1/4} \langle r \rangle^{-1/4} \partial \phi\|^2_{L^2(S_T)}\bigr),
 \nonumber \\
 \|\phi\|^2_{Z_2(T)}\!\! &=&\!\! \|\phi\|^2_{Z_1(T)} + \|\partial \phi\|^2_{Z_1(T)}.
 \nonumber
\end{eqnarray}
It is obvious that $Y_m(T) \subset Z_m(T)$. Moreover, we can easily see that
$C^2([0,T]\times {\mathbb R}^3) \cap X_m(T) \subset Y_m(T)$\ ($m = 1, 2$).

\bigskip

Now let us state our main result. Since we are concerned only with
radially symmetric solutions, we naturally suppose that $F(\partial \phi)$
takes the form
\begin{equation}
 F(\partial \phi) = a(\partial_t \phi)^2 + b|\nabla \phi|^2
\nonumber
\end{equation}
for some constants $a$, $b$.
After reading the proofs of our results, however, the readers will
be convinced that we can more generally consider the case where the
coefficients $a$ and $b$ depend on the unknown function $\phi$.


\begin{theorem}
 \label{ta1}
Suppose that $(f, g) \in H^2_{{\rm rad}}({\mathbb R}^3) \times H^1_{{\rm
 rad}}({\mathbb R}^3)$.
Then there exist positive constants
 $A_1, \varepsilon_0$ with the following property: if
$(f,g)$ and $T$ satisfy
\begin{equation}
 \|\nabla f\|_{H^1({\mathbb R}^3)} + \|g\|_{H^1({\mathbb R}^3)}=
  \varepsilon \leq \varepsilon_0
\quad \mbox{and}\quad T \leq \exp (A_1\varepsilon^{-1}),
\nonumber
\end{equation}
then the initial value problem {\rm (\ref{ea1})}, {\rm (\ref{ea2})} has a solution
 $\phi \in X_2(T) \cap Y_2(T)$.
Also, there exists a constant $M_1$ such that
\begin{equation}
 \|\phi\|_{E_2(T)}
 + \|\phi\|_{Y_2(T)}
 + \|\phi\|_{Z_2(T)}
\leq M_1 \varepsilon.
\nonumber
\end{equation}
Here the constants $A_1$, $\varepsilon_0$ and $M_1$ depend on $a$, $b$, and $h$.
\end{theorem}
We notice that, since $H^2({\mathbb R}^3) \subset L^{\infty}({\mathbb R}^3)$
and $H^1({\mathbb R}^3) \subset L^4({\mathbb R}^3)$,
the solution $\phi$ satisfies the equation (\ref{ea1})
as an equality in $C([0,T]; L^2({\mathbb R}^3))$.

\bigskip

We next consider the problem of well-posedness.
It is shown in Section 6.2 that the solution map $\Phi : (f, g) \mapsto \phi$
is well-defined on the set
\[
 D(\varepsilon_0) = \{(f,g)\,;\,
(f, g) \in H^2_{{\rm rad}}({\mathbb R}^3) \times
 H^1_{{\rm rad}}({\mathbb R}^3),
\|\nabla f\|_{H^1} + \|g\|_{H^1} < \varepsilon_0
\},
\]
where $\phi$ is a solution of (\ref{ea1}), (\ref{ea2}) with
$T = T_{\varepsilon} = \exp (A_1/\varepsilon)$, $\varepsilon = \|\nabla f\|_{H^1} + \|g\|_{H^1}$.
$D(\varepsilon_0)$ is a topological space with the relative topology
in $\dot{H}^1 \times L^2$.
Our result of well-posedness is as follows:
\begin{theorem}
 \label{ta2}
The solution map $\Phi: D(\varepsilon_0) \to X_1(T_{\varepsilon_0}) \cap Y_1(T_{\varepsilon_0})$ is well-defined
and Lipschitz continuous.
\end{theorem}

\bigskip

As mentioned above,
Hidano-Yokoyama \cite{HY} proved the almost global
existence of radial solutions for small initial data in $H^2 \times H^1$,
when $h = 0$ in (\ref{ea1}).
The key to this result is an effective use of the space-time $L^2$ estimate
\begin{eqnarray}
 \label{ea3}
 \lefteqn{\bigl(\log (2+T)\bigr)^{-1/2}\bigl(
\begin{array}[t]{l}
\|r^{-3/2+\mu}\langle r\rangle^{-\mu}\phi\|_{L^2((0,T)\times {\mathbb
R}^3)} \\
+  \|r^{-1/2+\mu}\langle r\rangle^{-\mu}\partial \phi\|_{L^2((0,T)\times
{\mathbb R}^3)}\bigr)
\rule{0pt}{3ex}
\end{array}
}
 \\
&\leq& C \Bigl(
\|\nabla f\|_{L^2} + \|g\|_{L^2} + \int_0^T\|\Box \phi (t,\cdot)\|_{L^2}dt
\Bigr),
\nonumber
\end{eqnarray}
where $r = |x|$, $\langle r \rangle = (1 + r^2)^{1/2}$ and $0 < \mu < 1/2$
(see also the Appendix of \cite{HY2}).
This type of inequalities is originally used in Keel-Smith-Sogge \cite{KSS},
and is very useful, because the right-hand side of
(\ref{ea3}) is the same as that of the energy inequality.
Combining this with the energy inequality, we can prove the almost
global existence for the semilinear case.
In order to adapt this approach in such a way that we can apply it to the present problem,
we need to prove a space-time $L^2$ estimate for the perturbed
wave equation with a variable coefficient differential operator
 (Theorem \ref{tb1}).
This generalization was made by Metcalfe-Sogge \cite{MS}.
They employed the method of Rodnianski \cite{Ster};
a space-time $L^2$ estimate is derived
from an energy-momentum tensor,
which produces several useful quantities by contracting it with suitable vector fields.
The vector field they used takes the form
\begin{equation}
 X = \sum_{a = 1}^n \frac{x^a}{\rho + r}\frac{\partial}{\partial x^a}.
\nonumber
\end{equation}
Since their $L^2_{t,x}$ estimate was for the study of the initial-boundary value problem in a domain exterior to obstacles,
we need to supplement their estimate with
an $L^2_{t,x}$ estimate near the spatial origin.
Following Metcalfe's suggestions \cite{M}, we use
\begin{equation}
 X = \sum_{a = 1}^n \left(\frac{r}{1 + r}\right)^{\kappa}
\frac{x^a}{r}\frac{\partial}{\partial x^a}
\nonumber
\end{equation}
for $0 < \kappa < 1$ and obtain a space-time $L^2$ estimate
on $(0,T)\times \{r < 1\}$. Combining this with the estimate of Metcalfe-Sogge \cite{MS},
we can obtain an $L^2_{t,x}$ estimate on $(0,T)\times {\mathbb R}^n$,
which plays a central role in our study.

As is usual with the quasilinear problems, we are faced with
difficulties caused by the loss of derivatives,
when solving the initial value problem (\ref{ea1}), (\ref{ea2}) by
successive approximations.
Because the loss of derivatives prevents us from showing
the convergence of the second derivatives,
it is far from trivial to prove that the solution thereby obtained
has the regularity stated in Theorem \ref{ta1}.
In order to prove the continuity in time
of the second derivatives,
we use the method of difference quotient.
We find in the proof that the space-time $L^2$ estimate is useful again.

We remark that the global existence fails in general, by results of
John \cite{J} and Sideris \cite{ST2}.
However, if we assume $F(\partial \phi) = 0$, we can expect the global
existence of the Cauchy problem (\ref{ea1}), (\ref{ea2}).
Lindblad \cite{L} proved the global existence of classical solutions for
small and smooth initial data, assuming radial symmetry.
When the radial symmetry is not imposed,
it is known by Alinhac \cite{A} that the global existence remains
true for small and smooth data. Furthermore, Lindblad \cite{L3} considered more general
quasilinear wave equations and proved the global existence result.
Note that all these results require considerably higher regularity and rapider decay at infinity for the initial data.
A significant improvement of this regularity assumption was made by Zhou-Lei \cite{ZL}, in which
they proved the global existence for radially symmetric, compactly supported
$H^2 \times H^1$-data.
Compared with that of \cite{ZL}, our analysis in this paper is not strong enough to prove the global existence,
but our result in Theorem \ref{ta1} does not require this support assumption.

This paper is organized as follows. In Section 2, we prove a
space-time $L^2$ estimate for a classical solution of the perturbed
wave equation.
Almost global existence for the initial value problem (\ref{ea1}),
(\ref{ea2}) is established from Section 3 through Section 5.
In Section 3, we define a sequence of successive approximation
to (\ref{ea1}), (\ref{ea2}) and prove its convergence to a weak solution
on some time interval.
It turns out to be a strong solution in Section 4,
and we further see that the solution actually exists almost globally in Section 5.
In Section 6, we discuss the problem of the continuity property of the solution map.
For this purpose, we need to confirm that the space-time $L^2$ estimate is
applicable to our low-regularity solutions. After that, we show that the
solution map is well-defined, and prove Theorem \ref{ta2}.

Throughout this paper, we use the summation convention that repeated
upper and lower indices are summed.
We use the Greek letters $\alpha$, $\beta$, $\dots$ when
the indices are from $0$ to $n$, and
the Latin letters $a$, $b$, $\dots$ when they run from $1$ to $n$.
We also use the geometric convention of raising and lowering indices
with respect to the Minkowski metric $(g_{\alpha\beta})={\rm diag}(-1,1,\dots,1)$.
As usual, $x^0 = t$ and $x = (x^1, x^2,\ldots , x^n)$ are the time and
the space variables respectively.
We remark $-\partial^0=\partial_0=\partial/\partial x^0$, $\partial^a=\partial_a
=\partial/\partial x^a$ for $a=1,2,\ldots,n$.
We set $r = |x|$ and $\langle r \rangle = \langle x \rangle = \sqrt{1 +
|x|^2}$. We let $|\cdot |$ denote the usual Euclidean norm,
namely $|x| = [(x^1)^2 + \cdots + (x^n)^2]^{1/2}$, $|\partial \phi|
= [(\partial_t \phi)^2 + (\partial_1 \phi)^2 + \cdots + (\partial_n \phi)^2]^{1/2}$ etc.


\section{Space-time $L^2$ estimate}

In this section, we consider the variable coefficient linear wave equation
\begin{equation}
 \partial_t^2 \phi - \Delta \phi + h^{\alpha \beta}(t,x)\partial_\alpha
  \partial_\beta \phi = F(t,x)
 \quad
\mbox{in}\,\,(0,T)\times{\mathbb R}^n,
\label{eb1}
\end{equation}
with initial data
\begin{equation}
 \phi(0,\cdot) = f \in H^1({\mathbb R}^n),\quad
\partial_t \phi (0,\cdot) = g \in L^2({\mathbb R}^n)
 \quad
\mbox{in}\,\,{\mathbb R}^n.
\label{eb2}
\end{equation}
Using the geometric multiplier method of Rodnianski in the Appendix of \cite{Ster},
we are going to prove a local-in-time space-time $L^2$ estimate
for a solution of (\ref{eb1}), (\ref{eb2}) when $n \geq 3$.

Let $h^{\alpha \beta} \in C^1([0,T] \times {\mathbb R}^n)
\enskip (\alpha, \beta = 0,1,\ldots ,n)$,
and suppose that they are uniformly bounded together with
their first derivatives.
In addition, we require that they satisfy the following conditions:
\begin{equation}
 \sum_{\alpha,\beta = 0}^n |h^{\alpha \beta}(t,x)| \leq \frac{1}{2}
 \quad
\mbox{for}\,\,(t,x) \in [0,T]\times {\mathbb R}^n,
\label{eb3}
\end{equation}
\begin{equation}
 h^{\alpha \beta}(t,x) = h^{\beta \alpha}(t,x)
 \quad
\mbox{for}\,\,
\begin{array}[t]{l}
 (t,x) \in [0,T]\times {\mathbb R}^n,\\
 \alpha,\beta = 0,1,2,\ldots, n.
\end{array}
\label{eb4}
\end{equation}


\begin{theorem}
\label{tb1}
Let $n \geq 3$.
Let $\phi \in C^2([0,T]\times {\mathbb R}^n)$
be a solution of the initial value problem {\rm (\ref{eb1})}, {\rm (\ref{eb2})}
satisfying the conditions {\rm (\ref{eb3})}, {\rm (\ref{eb4})}.
Suppose also that $\phi \in C^0([0,T];H^1({\mathbb R}^n)) \cap
 C^1([0,T];L^2({\mathbb R}^n))$ and
\[
 |\partial \phi||F| + \frac{|\phi||F|}{r^{1-2\mu}\langle r\rangle^{2\mu}}
 \in L^1((0,T)\times {\mathbb R}^n).
\]
 Then, for $0 < \mu < 1/2$, we have
\begin{eqnarray}
 \lefteqn{(1 + T)^{-2\mu}\Bigl(\|r^{-3/2 + \mu} \phi\|^2_{L^2((0,T) \times {\mathbb R}^n)}
+ \|r^{-1/2 + \mu} \partial \phi\|^2_{L^2((0,T) \times
{\mathbb R}^n)}\Bigr)} \nonumber \\
 \lefteqn{+\bigl(\log (2 + T)\bigr)^{-1}
\begin{array}[t]{l}
 \displaystyle
\Bigl(\|r^{-3/2 + \mu}\langle r\rangle^{-\mu} \phi\|^2_{L^2((0,T) \times {\mathbb R}^n)}
 \\
 \displaystyle
\quad+ \|r^{-1/2 + \mu}\langle r\rangle^{-\mu} \partial \phi\|^2_{L^2((0,T) \times
{\mathbb R}^n)}\Bigr)
\end{array}}
\label{eb5} \\
 &\leq& C_1(\|\nabla f\|^2_{L^2({\mathbb R}^n)} + \|g\|^2_{L^2({\mathbb
  R}^n)}) \nonumber \\
 & & + C_1\int_0^T\!\! \int_{{\mathbb R}^n} \left(
|\partial \phi||F| + \frac{|\phi||F|}{r^{1-2\mu}\langle r\rangle^{2\mu}} + |\partial h||\partial \phi|^2
 \right.
 \nonumber \\
 & & \qquad \quad \left.
+ \frac{|\partial h||\phi \partial \phi|}{r^{1-2\mu}\langle r\rangle^{2\mu}}
+ \frac{|h||\partial \phi|^2}{r^{1-2\mu}\langle r\rangle^{2\mu}}
+ \frac{|h||\phi \partial \phi|}{r^{2 - 2\mu}\langle r\rangle^{2\mu}} \right)dxdt.
 \nonumber
\end{eqnarray}
Here $C_1$ is a positive constant depending only on $n$ and $\mu$.
\end{theorem}

\bigskip

We begin by considering a weighted $L^2$
estimate of the solution of (\ref{eb1}), (\ref{eb2}) over the domain $(0,T)
\times \{x\,;\, x \in {\mathbb R}^n,\ |x| < 1\}$.


\begin{lemma}
 \label{lb1}
Let $n \geq 3$.
Let $\phi$ be the classical solution of {\rm (\ref{eb1})}, {\rm (\ref{eb2})} stated in Theorem {\rm \ref{tb1}}.
Then, for $0 < \mu < 1/2$, we have
\begin{eqnarray}
 \lefteqn{\int_0^T\!\! \int_{\{x \in {\mathbb R}^n ;\, |x| < 1\}} \Bigl(\frac{(\partial_t \phi)^2 +
  |\nabla \phi|^2}{r^{1-2\mu}} + \frac{\phi^2}{r^{3-2\mu}}\Bigr)dxdt}
 \label{eb6}\\
 &\leq& C(\|\nabla f\|^2_{L^2({\mathbb R}^n)} + \|g\|^2_{L^2({\mathbb
  R}^n)}) \nonumber \\
 & & + C\int_0^T\!\! \int_{{\mathbb R}^n} \left(
|\partial \phi||F| + \frac{|\phi||F|}{r^{1-2\mu}\langle r\rangle^{2\mu}} + |\partial h||\partial \phi|^2
 \right.
 \nonumber \\
 & & \qquad \quad \left.
+ \frac{|\partial h||\phi \partial \phi|}{r^{1-2\mu}\langle r\rangle^{2\mu}}
+ \frac{|h||\partial \phi|^2}{r^{1-2\mu}\langle r\rangle^{2\mu}}
+ \frac{|h||\phi \partial \phi|}{r^{2 - 2\mu}\langle r\rangle^{2\mu}} \right)dxdt.
 \nonumber
\end{eqnarray}
Here $C$ is a positive constant depending only on $n$ and $\mu$.
\end{lemma}

Proof.
Following \cite{MS}, we define the energy-momentum tensor
\begin{equation}
Q_{\alpha\beta}[\phi] =
\begin{array}[t]{l}
\displaystyle
\partial_\alpha\phi\partial_\beta\phi
-\frac12g_{\alpha\beta}\partial^\gamma\phi\partial_\gamma\phi \\
\displaystyle
- g_{\alpha\gamma}h^{\gamma\delta}\partial_\delta\phi\partial_\beta\phi
+
\frac12g_{\alpha\beta}h^{\gamma\delta}\partial_\gamma\phi\partial_\delta
\phi
\end{array}
\label{eb7}
\end{equation}
for $\alpha, \beta = 0,1,\ldots ,n$,
where $(g_{\alpha\beta})={\rm diag}(-1,1,\dots,1)$.
For a function $f(r)$ suitably chosen later,
we set
\begin{equation}
X^0=0,\,\,\,
X^a = \frac{f(r)x^a}{r}
\quad \mbox{for}\,\,a=1,2,\dots,n.
\label{eb8}
\end{equation}
Contracting the energy-momentum tensor $Q[\phi]$
with the vector field $X$, we define
\begin{equation}
P_\alpha[\phi,X]
=
Q_{\alpha\beta}[\phi]X^\beta
\quad \mbox{for}\,\,\alpha = 0,1,\ldots, n.
\label{eb9}
\end{equation}
We also set
\begin{equation}
 \pi_{\alpha\beta} = \frac{1}{2}(\partial_{\alpha}X_{\beta} + \partial_{\beta}X_{\alpha})
\quad
\mbox{for}\,\,\alpha,\beta = 0,1,\ldots, n.
\nonumber
\end{equation}
Then we have
\begin{eqnarray}
 \lefteqn{\partial^{\alpha}P_{\alpha}[\phi,X]}
\nonumber \\
&=& f'(r)(\partial_r \phi)^2 + \frac{f(r)}{r}|\not{\!\!\nabla}\phi|^2 -
 \frac{1}{2}\mathrm{tr}\pi\,\partial^{\gamma}\phi\partial_{\gamma}\phi
 - f(r)(\partial_r \phi)F
\nonumber \\
& & - f(r)\partial_{\gamma}h^{\gamma\delta}\partial_{\delta}\phi\partial_r\phi
+\frac{1}{2}f(r)\partial_rh^{\gamma\delta}\partial_{\gamma}\phi\partial_{\delta}\phi
-\frac{x^a}{r}f'(r)h^{a\delta}\partial_{\delta}\phi\partial_r\phi
 \nonumber \\
& & +\frac{x_a f(r)}{r^2}h^{a\delta}\partial_{\delta}\phi\partial_r\phi
-\frac{f(r)}{r}h^{a\delta}\partial_{\delta}\phi\partial_a\phi
+\frac{1}{2}(\mathrm{tr}\pi)h^{\gamma\delta}\partial_{\gamma}\phi\partial_{\delta}\phi,
\nonumber
\end{eqnarray}
where
\begin{equation}
 \mathrm{tr}\pi = f'(r) + (n - 1)\frac{f(r)}{r}
\nonumber
\end{equation}
(see \cite{MS}, p.199).
Here $\not{\!\!\nabla}\phi$ denotes the angular portion
of the spatial gradient $\nabla \phi$.
Furthermore, we introduce the modified momentum density
\begin{eqnarray}
 {\bar P}_{\alpha}[\phi,X] &=&
P_{\alpha}[\phi,X]
+ \frac{n - 1}{2}\left(\frac{f(r)}{r}\right)\phi\partial_{\alpha}\phi
 \label{eb10} \\
& & - \frac{n - 1}{4}\partial_{\alpha}\left(\frac{f(r)}{r}\right)\phi^2
- \frac{n - 1}{2}\left(\frac{f(r)}{r}\right)g_{\alpha\gamma}h^{\gamma\beta}\phi\partial_{\beta}\phi
 \nonumber
\end{eqnarray}
for $\alpha = 0,1,\ldots ,n$.
Then we obtain
\begin{eqnarray}
 \partial^{\alpha} {\bar P}_{\alpha}[\phi,X] &=&
 f'(r)(\partial_r \phi)^2 + \frac{f(r)}{r}|\not{\!\!\nabla}\phi|^2
 - \frac{1}{2}f'(r)\partial^{\gamma}\phi\partial_{\gamma}\phi
\nonumber  \\
 & & - \frac{n-1}{4}\Delta \left(\frac{f(r)}{r}\right) \phi^2 + {\bar R},
\label{eb11}
\end{eqnarray}
where
\begin{eqnarray}
 {\bar R} &=& -f(r)(\partial_r \phi)F -
  \frac{n-1}{2}\left(\frac{f(r)}{r}\right)\phi F \label{eb12} \\
 & & - f(r)(\partial_{\gamma}h^{\gamma\delta})\partial_{\delta}\phi \left(
\partial_r \phi + \frac{n-1}{2}\frac{\phi}{r}\right) \nonumber \\
 & & + \frac{1}{2}f(r)(\partial_r h^{\gamma\delta})\partial_{\gamma}\phi
  \partial_{\delta}\phi \nonumber \\
 & & - \frac{x_a f'(r)}{r}h^{a\delta}\partial_{\delta}\phi\left(\partial_r \phi +
  \frac{n-1}{2} \frac{\phi}{r}\right) \nonumber \\
 & & + \frac{x_a f(r)}{r^2}h^{a\delta}\partial_{\delta}\phi\left(\partial_r \phi +
  \frac{n-1}{2} \frac{\phi}{r}\right) \nonumber \\
 & & - \frac{f(r)}{r}h^{a\delta}\partial_{\delta}\phi \partial_a \phi
 + \frac{1}{2}f'(r) h^{\gamma\delta}\partial_{\gamma}\phi \partial_{\delta} \phi \nonumber
\end{eqnarray}
(see \cite{MS}, pp.199--200).
The identity (\ref{eb11}) further leads to
\begin{eqnarray}
 \lefteqn{\partial^{\alpha} {\bar P}_{\alpha}[\phi,X] =
 \frac{1}{2}f'(r)(|\nabla \phi|^2+|\partial_t \phi|^2)}
 \label{eb13} \\
 & & + \left(\frac{f(r)}{r} - f'(r)\right)|\not{\!\!\nabla}\phi|^2
- \frac{n-1}{4}\Delta \left(\frac{f(r)}{r}\right) \phi^2 + {\bar R}.
\nonumber
\end{eqnarray}

Now we set
\begin{equation}
 f(r) = \left(\frac{r}{1 + r}\right)^{\kappa}
 \label{eb14}
\end{equation}
for $0 < \kappa < 1$. We first show that
\begin{equation}
 \frac{f(r)}{r} - f'(r) \geq
 (1 - \kappa) \frac{r^{\kappa - 1}}{(1 + r)^{\kappa}}
\label{eb15}
\end{equation}
and
\begin{equation}
 \Delta \left(\frac{f(r)}{r}\right) \leq
 -\frac{\kappa (1 - \kappa)}{r^{3 - \kappa}(1 + r)^{2 + \kappa}}.
\label{eb16}
\end{equation}
In fact, (\ref{eb15}) is an immediate result of
\[
 \frac{f(r)}{r} - f'(r) = \frac{r^{\kappa - 1}}{(1 + r)^{\kappa}}\left(
1 - \frac{\kappa}{1 + r}\right).
\]
In order to compute $\Delta (f(r)/r)$, we note that
\begin{eqnarray*}
 \Delta \left(\frac{f(r)}{r}\right)
 &=& r^{1-n}\partial_r\left(r^{n-1}\partial_r \frac{f(r)}{r}\right) \\
 &=& r^{1-n}\partial_r\left(r^{n-2}\left(f'(r)-\frac{f(r)}{r}\right)\right).
\end{eqnarray*}
Using this identity, we have
\begin{eqnarray*}
 \Delta \left(\frac{f(r)}{r}\right)
 &=& r^{1-n}\partial_r\left(\frac{r^{n-3+\kappa}}{(1+r)^{\kappa}}
 \left(\frac{\kappa}{1+r}-1\right)\right) \\
 &=& \left(\frac{(n-3+\kappa)r^{-3+\kappa}}{(1+r)^{\kappa}}
-\frac{\kappa r^{-2+\kappa}}{(1+r)^{\kappa+1}}\right)
 \left(\frac{\kappa}{1+r}-1\right) \\
 & &- \frac{\kappa r^{-2+\kappa}}{(1+r)^{\kappa + 2}},
\end{eqnarray*}
which gives (\ref{eb16}), because
\begin{eqnarray*}
 \Delta \left(\frac{f(r)}{r}\right)\!\! &\leq&\!\!
\left(\frac{\kappa r^{-3+\kappa}}{(1+r)^{\kappa}}
-\frac{\kappa r^{-2+\kappa}}{(1+r)^{\kappa+1}}\right)
 \left(\frac{\kappa}{1+r}-1\right)
\\
 &\leq & -\frac{\kappa (1 - \kappa)}{r^{3 - \kappa}(1 + r)^{2 + \kappa}}.
 \nonumber
\end{eqnarray*}

If we integrate (\ref{eb13}) over $[0,T] \times {\mathbb R}^n$ and apply
the divergence theorem, we get
\begin{eqnarray}
\lefteqn{%
 \int_0^T\!\! \int_{{\mathbb R}^n}\left[
 \frac{1}{2}f'(r)(|\nabla \phi|^2+|\partial_t \phi|^2)
  + \left(\frac{f(r)}{r} - f'(r)\right)|\not{\!\!\nabla}\phi|^2 \right.}
 \nonumber \\
 & & \qquad \left.  - \frac{n-1}{4}\Delta \left(\frac{f(r)}{r}\right) \phi^2 \right]dxdt
 \label{eb17} \\
 &=& - \int_{{\mathbb R}^n} {\bar P}_0[\phi,X](T,x)dx
 + \int_{{\mathbb R}^n} {\bar P}_0[\phi,X](0,x)dx
- \int_0^T\!\! \int_{{\mathbb R}^n} {\bar R}\, dxdt.
 \nonumber
\end{eqnarray}
To obtain (\ref{eb17}),
we should note that there exists a sequence $\{R_k\}$ so that
$R_k \to \infty$ and
\begin{equation}
 \lim_{k \to \infty} \int_0^T\!\! \int_{|x| = R_k} {\bar P}_a[\phi,X](t,x)\cdot
  \frac{x^a}{r}\,dSdt = 0.
 \label{eb18}
\end{equation}
Indeed, the existence of such a sequence can be easily deduced from
${\bar P}_a[\phi,X]\cdot x^a/r \in L^1([0,T]\times {\mathbb R}^n)$.
To see this integrability,
we recall (\ref{eb7})--(\ref{eb10}) and obtain
\begin{eqnarray*}
 {\bar P}_a[\phi,X]\cdot \frac{x^a}{r}
 &=& P_a[\phi,X]\cdot \frac{x^a}{r} + \frac{n-1}{2}\frac{f(r)}{r}\phi
 \partial_r \phi \\
 & & - \frac{n-1}{4} \partial_r \left(\frac{f(r)}{r}\right) \phi^2
 - \frac{n-1}{2}\frac{x_af(r)}{r^2}h^{a\beta} \phi \partial_{\beta} \phi \\
 &=& f(r)(\partial_r \phi)^2 - \frac{1}{2}f(r)\partial^{\gamma} \phi \partial_{\gamma} \phi
 - \frac{x_a f(r)}{r}h^{a\delta}\partial_{\delta} \phi \partial_r \phi \\
 & & + \frac{1}{2}f(r)h^{\gamma\delta}\partial_{\gamma}\phi \partial_{\delta} \phi
+ \frac{n-1}{2}\frac{f(r)}{r}\phi
 \partial_r \phi \\
 & & - \frac{n-1}{4} \partial_r \left(\frac{f(r)}{r}\right) \phi^2
 - \frac{n-1}{2}\frac{x_af(r)}{r^2}h^{a\beta} \phi \partial_{\beta} \phi.
\end{eqnarray*}
This gives the estimate
\[
 \left|{\bar P}_a[\phi,X]\cdot \frac{x^a}{r}\right| \leq C\left(|\partial \phi|^2 +
 \frac{\phi^2}{r^2}\right),
\]
because $|f(r)| \leq 1$, $|\partial_r \bigl(f(r)/r\bigr)| \leq Cr^{-2}$
by (\ref{eb14}). Hence we conclude ${\bar P}_a[\phi,X] \cdot x^a/r \in
L^1([0,T]\times {\mathbb R}^n)$ by Hardy's inequality.

Going back to (\ref{eb17}) and using (\ref{eb15}), (\ref{eb16}), we have
\begin{eqnarray}
 \lefteqn{%
 \int_0^T \int_{{\mathbb R}^n}\left[
\frac{(\partial_t \phi)^2 + |\nabla \phi|^2}{r^{1-\kappa}(1 + r)^{1 + \kappa}}
+\frac{|\not{\!\!\nabla} \phi|^2}{r^{1-\kappa}(1 + r)^{\kappa}}
+\frac{\phi^2}{r^{3-\kappa}(1 + r)^{2 + \kappa}}
\right]dxdt
} \nonumber \\
&\leq & C\left(\|{\bar P}_0[\phi,X](T,\cdot)\|_{L^1({\mathbb R}^n)}
+\|{\bar P}_0[\phi,X](0,\cdot)\|_{L^1({\mathbb R}^n)}
\right) \label{eb19} \\
& & + C \int_0^T \int_{{\mathbb R}^n} |{\bar R}|\,dxdt. \nonumber
\end{eqnarray}
So it remains to bound the right-hand side of (\ref{eb19}).
By (\ref{eb7})--(\ref{eb10}), we get
\begin{eqnarray}
 {\bar P}_0[\phi,X] &=& P_0[\phi,X]
 + \frac{n-1}{2}\frac{f(r)}{r}\phi \partial_t\phi
+  \frac{n-1}{2}\frac{f(r)}{r}h^{0\beta}\phi \partial_{\beta}\phi
\nonumber \\
&=& f(r) \partial_t \phi \partial_r \phi + f(r)h^{0\delta} \partial_{\delta} \phi \partial_r \phi
\nonumber \\
& & +  \frac{n-1}{2}\frac{f(r)}{r}\phi \partial_t\phi
+  \frac{n-1}{2}\frac{f(r)}{r}h^{0\beta}\phi \partial_{\beta}\phi.
\nonumber
\end{eqnarray}
Since $|f(r)| \leq 1$ and $|h^{\alpha \beta}| \leq 1/2$,
we have
\begin{eqnarray}
 \label{eb20}
 \lefteqn{\|{\bar P}_0[\phi,X](T,\cdot)\|_{L^1({\mathbb R}^n)}}\\
 &\leq& C \left(
\|\partial \phi(T,\cdot)|^2_{L^2({\mathbb R}^n)} +
\left\|\frac{\phi(T,\cdot)}{r}\right\|^2_{L^2({\mathbb R}^n)}\right)
\nonumber \\
 &\leq&  C \|\partial \phi(T,\cdot)\|^2_{L^2({\mathbb R}^n)}.
\nonumber
\end{eqnarray}
Similarly,
\begin{equation}
 \|{\bar P}_0[\phi,X](0,\cdot)\|_{L^1({\mathbb R}^n)} \leq
C \|\partial \phi(0,\cdot)\|^2_{L^2({\mathbb R}^n)}.
\label{eb21}
\end{equation}
In addition, the standard energy inequality yields
\begin{eqnarray}
 \|\partial \phi(T,\cdot)\|^2_{L^2({\mathbb R}^n)}
& \leq & C \|\partial \phi(0,\cdot)\|^2_{L^2({\mathbb R}^n)}
 \label{eb22} \\
& & + C \int_0^T\!\! \int_{{\mathbb R}^n} (|\partial_t \phi F| + |\partial
 h||\partial \phi|^2)\,dxdt.
\nonumber
\end{eqnarray}
Lastly, we see from (\ref{eb12}) and (\ref{eb14}) that
\begin{eqnarray}
 |{\bar R}| &\leq &
 C\left(|\partial \phi||F| + \frac{|\phi||F|}{r^{1-\kappa}(1 + r)^{\kappa}} + |\partial
   h||\partial \phi|^2 \right. \label{eb23} \\
 & & \left. + \frac{|\partial h||\phi \partial \phi|}{r^{1-\kappa}(1 + r)^{\kappa}} +
  \frac{|h||\partial \phi|^2}{r^{1-\kappa}(1 + r)^{\kappa}}
 + \frac{|h||\phi \partial \phi|}{r^{2-\kappa}(1 + r)^{\kappa}}\right).
 \nonumber
\end{eqnarray}
Combining (\ref{eb19})--(\ref{eb23}),
we conclude that
\begin{eqnarray}
 \lefteqn{%
 \int_0^T \int_{\{x \in {\mathbb R}^n ;\, |x|<1\}} \left(
\frac{|\partial \phi|^2}{r^{1-\kappa}} + \frac{\phi^2}{r^{3 - \kappa}}
\right)\,dxdt}
\nonumber \\
 &\leq & C \|\partial \phi(0,\cdot)\|^2_{L^2({\mathbb R}^n)}
\nonumber \\
 & & + C \int_0^T \int_{{\mathbb R}^n} \left(
|\partial \phi||F| + \frac{|\phi||F|}{r^{1-\kappa}(1 + r)^{\kappa}}
+ |\partial h||\partial \phi|^2 \right.
\nonumber \\
& & \left. + \frac{|\partial h||\phi \partial \phi|}{r^{1-\kappa}(1 + r)^{\kappa}} +
 \frac{|h||\partial \phi|^2}{r^{1-\kappa}(1+r)^{\kappa}} + \frac{|h||\phi \partial \phi|}{r^{2-\kappa}(1+r)^{\kappa}}
\right)\,dxdt. \nonumber
\end{eqnarray}
This gives (\ref{eb6}), if we set $\kappa = 2\mu
\enskip (0 < \mu < 1/2)$.
\hfill $\Box$


\begin{lemma}
 \label{lb2}
Let $n \geq 3$.
Let $\phi$ be the classical solution of {\rm (\ref{eb1})}, {\rm (\ref{eb2})} stated in Theorem {\rm \ref{tb1}}.
Then, for $T > 1$ and $0 < \mu < 1/2$, we have
\begin{eqnarray}
 \lefteqn{T^{-2\mu}\!\!\!\int_0^T\!\!\! \int_{\{x \in {\mathbb R}^n ;\, 1 < r < T\}}\! \Bigl(\frac{(\partial_t \phi)^2 +
  |\nabla \phi|^2}{r^{1-2\mu}} + \frac{\phi^2}{r^{3-2\mu}}\Bigr)dxdt}
 \label{eb24}\\
 \lefteqn{+ \bigl(\log (2+T)\bigr)^{-1}\!\!\!\int_0^T\!\!\!
\int_{\{x \in {\mathbb R}^n ;\, 1 < r < T\}}\! \Bigl(\frac{(\partial_t \phi)^2 +
  |\nabla \phi|^2}{r} + \frac{\phi^2}{r^3}\Bigr)dxdt}
 \nonumber \\
 &\leq& C(\|\nabla f\|^2_{L^2({\mathbb R}^n)} + \|g\|^2_{L^2({\mathbb
  R}^n)}) \nonumber \\
 & & + C\int_0^T \int_{{\mathbb R}^n} \left(
|\partial \phi||F| + \frac{|\phi||F|}{\langle r\rangle} + |\partial h||\partial \phi|^2
 \right.
 \nonumber \\
 & & \qquad \quad \left.
+ \frac{|\partial h||\phi \partial \phi|}{\langle r\rangle}
+ \frac{|h||\partial \phi|^2}{\langle r\rangle}
+ \frac{|h||\phi \partial \phi|}{\langle r\rangle^2} \right)dxdt.
 \nonumber
\end{eqnarray}
Here $C$ is a positive constant depending only on $n$ and $\mu$.
\end{lemma}

Proof. Though this estimate is essentially proven in \cite{MS}, we show the
proof for the sake of completeness.
Just as in the proof of the last lemma, we use the geometric
multiplier method. Instead of (\ref{eb14}), we choose
\begin{equation}
 f(r) = \frac{r}{\rho + r}
\label{eb25}
\end{equation}
(see \cite{MS}, p.197). If $\rho = 2^k$ and $2^{k - 1} \leq r \leq 2^k$
$(k = 1,2,\ldots)$,
we can easily see that
\begin{eqnarray}
\label{eb26}
 f'(r) &=& \frac{\rho}{(\rho + r)^2} \geq \frac{1}{2(\rho + r)}, \\
\label{eb27}
 \frac{f(r)}{r} - f'(r) &=& \frac{r}{(\rho + r)^2} \geq \frac{1}{3(\rho + r)},
\end{eqnarray}
and
\begin{eqnarray}
 -\Delta \left(\frac{f(r)}{r}\right) &=& \frac{(n-3)r + (n-1)\rho}{r(\rho + r)^3}
\label{eb28} \\
 &\geq & \frac{(n-1)\rho}{r(\rho + r)^3} \geq \frac{n-1}{(\rho + r)^3}.
\nonumber
\end{eqnarray}
Combining (\ref{eb17}) and (\ref{eb26})--(\ref{eb28}), we get
\begin{eqnarray}
 \lefteqn{%
 \int_0^T \int_{2^{k-1} < r < 2^k} \left(
 \frac{(\partial_t \phi)^2 + |\nabla \phi|^2}{2^k + r}
 + \frac{\phi^2}{(2^k + r)^3}
\right)\,dxdt} \label{eb29} \\
&\leq & 4\left( \|{\bar P}_0[\phi,X](T,\cdot)\|_{L^1({\mathbb R}^n)}
 + \|{\bar P}_0[\phi,X](0,\cdot)\|_{L^1({\mathbb R}^n)}\right)
\nonumber \\
& & + 4 \int_0^T \int_{{\mathbb R}^n} |{\bar R}|\,dxdt
\quad \mbox{for}\,\, k = 1,2,\ldots.
\nonumber
\end{eqnarray}

We note that the $L^1({\mathbb R}^n)$-norm of ${\bar P}_0[\phi,X](T,x)$ and that of ${\bar P}_0[\phi,X](0,x)$ have a bound similar to
({\ref{eb20}}) and ({\ref{eb21}}), respectively,
with the constant $C$ independent of $k$. Concerning ${\bar R}$, we see from (\ref{eb12}) that
\begin{eqnarray*}
 |{\bar R}| &\leq & C \left(
 |\partial \phi||F| + \frac{|\phi||F|}{1 + r} + |\partial h||\partial \phi|^2 \right.
\\
 & & \left. + \frac{|\partial h||\phi \partial \phi|}{1 + r} + \frac{|h||\partial \phi|^2}{1 + r}
 + \frac{|h||\phi \partial \phi|}{(1 + r)^2}
\right)
\end{eqnarray*}
for a positive constant $C$ which is independent of $k$.
Therefore, it follows from (\ref{eb29}) that
\begin{eqnarray}
 \lefteqn{%
 \int_0^T \int_{2^{k-1} < r < 2^k} \left(
 \frac{(\partial_t \phi)^2 + |\nabla \phi|^2}{r}
 + \frac{\phi^2}{r^3}
\right)\,dxdt} \label{eb30} \\
&\leq & C \|\partial \phi(0,\cdot)\|^2_{L^2({\mathbb R}^n)}
 \nonumber \\
& & + C \int_0^T \int_{{\mathbb R}^n}
 \left(
 |\partial \phi||F| + \frac{|\phi||F|}{\langle r \rangle} + |\partial h||\partial \phi|^2 \right.
 \nonumber \\
 & &
 \qquad
 \left. + \frac{|\partial h||\phi \partial \phi|}{\langle r \rangle} + \frac{|h||\partial \phi|^2}{\langle r \rangle}
 + \frac{|h||\phi \partial \phi|}{\langle r \rangle^2}
\right)\,dxdt. \nonumber
\end{eqnarray}

Now let $N$ be the smallest integer not smaller than $\log_2 T$. If $0<\mu<1/2$,
we have
\begin{eqnarray}
  \lefteqn{%
T^{-2\mu}\!\!\! \int_0^T\!\!\! \int_{1 < r < T}\! \left(
 \frac{|\partial \phi|^2}{r^{1-2\mu}}
 + \frac{\phi^2}{r^{3-2\mu}}\right)\,dxdt} \label{eb31} \\
 & &+ \bigl(\log (2 + T)\bigr)^{-1}\!\!\!
\int_0^T\!\!\! \int_{1 < r < T}\! \left(
 \frac{|\partial \phi|^2}{r}
 + \frac{\phi^2}{r^3}\right)\,dxdt
 \nonumber \\
&\leq &T^{-2\mu}\! \sum_{k=1}^{N}\!\! \int_0^T\!\!\! \int_{2^{k-1} < r < 2^k}\!\! \left(
 \frac{|\partial \phi|^2}{r^{1-2\mu}}
 + \frac{\phi^2}{r^{3-2\mu}}\right)\,dxdt \nonumber \\
 & &+ \bigl(\log (2 + T)\bigr)^{-1}\! \sum_{k=1}^{N}\!\!
\int_0^T\!\!\! \int_{2^{k-1} < r < 2^k}\!\! \left(
 \frac{|\partial \phi|^2}{r}
 + \frac{\phi^2}{r^3}\right)\,dxdt
 \nonumber \\
&\leq &T^{-2\mu}\! \sum_{k=1}^{N}(2^k)^{2\mu}\!\! \int_0^T\!\!\! \int_{2^{k-1} < r < 2^k}\!\! \left(
 \frac{|\partial \phi|^2}{r}
 + \frac{\phi^2}{r^3}\right)\,dxdt \nonumber \\
 & &+ \bigl(\log (2 + T)\bigr)^{-1}\! \sum_{k=1}^{N}\!\!
\int_0^T\!\!\! \int_{2^{k-1} < r < 2^k}\!\! \left(
 \frac{|\partial \phi|^2}{r}
 + \frac{\phi^2}{r^3}\right)\,dxdt
 \nonumber \\
 &\leq & C \sup_{k \geq 1} \int_0^T\!\!\! \int_{2^{k-1} < r < 2^k}\!\! \left(
 \frac{|\partial \phi|^2}{r}
 + \frac{\phi^2}{r^3}\right)\,dxdt.
\nonumber
\end{eqnarray}
By (\ref{eb30}) and (\ref{eb31}), we obtain (\ref{eb24}).
\hfill $\Box$

\bigskip

Proof of Theorem \ref{tb1}.
In view of Lemma \ref{lb1} and Lemma \ref{lb2},
it remains to prove the weighted $L^2$ estimate of the solution over
$(0,T) \times \{x\,;\, x \in {\mathbb R}^n,\ |x| > T\}$. If $0 < \mu < 1/2$, we have
\begin{eqnarray*}
  \lefteqn{%
(1 + T)^{-2\mu}\!\!\! \int_0^T\!\!\! \int_{r > T}\! \left(
 \frac{|\partial \phi|^2}{r^{1-2\mu}}
 + \frac{\phi^2}{r^{3-2\mu}}\right)\,dxdt} \\
 & &+ \bigl(\log (2 + T)\bigr)^{-1}\!\!\!
\int_0^T\!\!\! \int_{r > T}\! \left(
 \frac{|\partial \phi|^2}{r^{1-2\mu}\langle r \rangle^{2\mu}}
 + \frac{\phi^2}{r^{3-2\mu}\langle r \rangle^{2\mu}}\right)\,dxdt
 \nonumber \\
&\leq &\frac{C}{T}\! \int_0^T\!\!\! \int_{r > T}\!\! \left(
 |\partial \phi|^2
 + \frac{\phi^2}{r^2}\right)\,dxdt \nonumber \\
&\leq & C \sup_{0 < t < T} \int_{{\mathbb R}^n}\!\! \left(
 |\partial \phi|^2
 + \frac{\phi^2}{r^2}\right)\,dx. \nonumber
\end{eqnarray*}
Using Hardy's inequality and the standard energy inequality,
we see that the last quantity is bounded by the right-hand side of (\ref{eb5}).
\hfill $\Box$


\section{Local Existence}

\subsection{Preliminaries}

This section is concerned with the initial value problem for
the quasilinear wave equation
(\ref{ea1}).
We ask for a solution of (\ref{ea1}) with initial values (\ref{ea2}).
For this purpose,
we define a sequence of functions as follows.

Firstly, we choose a suitable radially symmetric function
$\rho \in C_0^{\infty}({\mathbb R}^3)$ and set
\begin{equation}
 \rho_j (x) = j^3 \rho (jx)
\quad \mbox{for}\,\, j = 1,2,3,\ldots
\label{ed7}
\end{equation}
so that
\begin{eqnarray}
 & & \rho_j \in C_0^{\infty}({\mathbb R}^3),\quad
 \mathrm{supp}\, \rho_j \subset
 \left\{x\,;\, x \in {\mathbb R}^3,\ |x| < \frac{1}{j}\right\},
\nonumber \\
 & & \int_{{\mathbb R}^3} \rho_j(x)\,dx = 1,\quad
 \rho_j(x) \geq 0 \quad \mbox{for}\,\, j=1,2,3,\ldots
\nonumber
\end{eqnarray}

\bigskip

Secondly, we set
\begin{equation}
 f_k(x) = \rho_{2^k} * f(x),\quad
g_k(x) = \rho_{2^k} * g(x)\quad
\label{ed11}
\end{equation}
for $k = 0, 1, 2,\ldots$
As is well known,
\begin{equation}
 f_k,\ g_k \in C^{\infty}({\mathbb R}^3)\quad
\mbox{for}\,\, k = 0,1,2,\ldots,
\nonumber
\end{equation}
\begin{equation}
 \|f_k - f\|_{H^2} \rightarrow 0,\quad
 \|g_k - g\|_{H^1} \rightarrow 0 \quad
\mbox{as}\,\, k \to \infty.
\nonumber
\end{equation}
Note that if $f$ and $g$ are radially symmetric, then $f_k$ and $g_k$
are also radially symmetric.
We also have
\begin{equation}
 \sum_{k = 1}^{\infty} \bigl( \|\nabla f_k
 - \nabla f_{k-1}\|_{L^2({\mathbb R}^3)}
+ \|g_k
 - g_{k-1}\|_{L^2({\mathbb R}^3)}\bigr)
< \infty.
\label{ed12}
\end{equation}
Indeed, we can easily check this by using
\begin{equation}
\|\rho_{2^k} * \varphi - \varphi\|_{L^2} \leq C2^{-k}\|\varphi\|_{H^1}.
\nonumber
\end{equation}

Finally, we let $\phi_{-1} \equiv 0$ and define $\phi_k$ $(k =
0,1,2,\ldots)$ recursively by solving
\begin{eqnarray}
 & &\partial_t^2 \phi_k - \Delta \phi_k + h(\phi_{k-1})
  \Delta \phi_k =  F(\partial \phi_{k-1})
\quad \mbox{in}\,\, S_T,
 \label{ed14} \\
 & &\phi_k(0,\cdot) = f_k,\quad
 \partial_t \phi_k(0,\cdot) = g_k \quad
 \mbox{in}\,\, {\mathbb R}^3.
 \label{ed15}
\end{eqnarray}
In order to ensure that the sequence $\{\phi_k\}$ is well-defined,
we will have to assume that $\varepsilon = \|\nabla f\|_{H^1} + \|g\|_{H^1}$ is small enough
and $T \leq \exp (A_1 \varepsilon^{-1})$ for some positive constant $A_1$
(see Lemma \ref{ld1}).
Applying a standard existence, uniqueness and regularity theorem to
(\ref{ed14})--(\ref{ed15}), we will see that, for all $k = 0,1,2,\ldots$,
$\phi_k$ is certainly defined,
radially symmetric for all time and satisfies
\begin{equation}
\phi_k \in C^{\infty}(\overline{S_T}) \cap X_2(T) \bigl(\subset
 Y_2(T) \bigr).
  \label{ed16}
\end{equation}

In order to estimate the sequence,
we will frequently use the Sobolev-type inequalities below.


\begin{lemma}
 \label{ld2}
Assume that $\phi \in \cap_{j=0}^2 C^j([0,T]; H^{2-j}({\mathbb R}^3))$.
Then we have
\begin{equation}
 |\phi(t,x)| \leq C_{\rm S} \|\phi\|_{E_2(T)} \quad
 \label{ed21}
\end{equation}
for $(t,x) \in \overline{S_T}$.
If in addition $\phi(t,\cdot)$ is radially symmetric, we also have
\begin{eqnarray}
 & & |\phi(t, x)| \leq C_{\rm S}\, r^{-1/2} \|\nabla \phi(t,\cdot)\|_{L^2({\mathbb R}^3)},
 \label{ed23} \\
 & & |\partial \phi(t,x)| \leq C_{\rm S}\, r^{-1/2}\langle r \rangle^{-1/2}
  \|\phi\|_{E_2(T)}
 \label{ed22}
\end{eqnarray}
for $(t,x) \in [0,T]\times ({\mathbb R}^3\setminus \{0\})$.
Here $C_{\rm S}$ is a positive constant.
\end{lemma}

Proof. The inequaliy (\ref{ed21}) is a result of the Sobolev embeddings
$W^{1,6}({\mathbb R}^3) \subset L^{\infty}({\mathbb R}^3)$
and ${\dot H^1}({\mathbb R}^3) \subset L^6({\mathbb R}^3)$.
To obtain (\ref{ed23}) and (\ref{ed22}), we use the Sobolev-type inequalities
\begin{eqnarray}
& & r^{1/2}|v(x)| \leq C \sum_{|\alpha| \leq 1} \|\nabla \Omega^{\alpha}
  v\|_{L^2({\mathbb R}^3)},
 \label{ed23e} \\
& & r|v(x)| \leq C \sum_{|\alpha| \leq 2} \|\Omega^{\alpha}
  v\|_{L^2({\mathbb R}^3)} + C \sum_{|\alpha| \leq 1} \|\nabla \Omega^{\alpha}
  v\|_{L^2({\mathbb R}^3)}.
 \label{ed24}
\end{eqnarray}
See Lemma 4.2 of \cite{KS} and Lemma 3.3 of \cite{ST}.
Here $\Omega = (\Omega_{12}, \Omega_{23}, \Omega_{31})$, $\Omega_{ab} =
x_a \partial_b - x_b\partial_a$ $(a, b = 1,2,3)$. Note that $\Omega v =
0$ if $v$ is radially symmetric.
Hence, (\ref{ed23}) immediately follows from (\ref{ed23e}).
If we use (\ref{ed23e}) for $r<1$, and
(\ref{ed24}) for $r>1$, we get
\[
 |v(x)| \leq C r^{-1/2}\langle r \rangle^{-1/2}
  \|v\|_{H^1({\mathbb R}^3)}
\]
for a radially symmetric function $v$.
Applying this inequality
to $\partial_t \phi$, $\partial_r \phi$
and noting $\|\partial_r \phi(t,\cdot)\|_{H^1} \leq
C\|\nabla \phi(t,\cdot)\|_{H^1}$,
we have (\ref{ed22}).
\hfill $\Box$

\bigskip

If we let $n = 3$ and $\mu = 1/4$ in Theorem \ref{tb1},
(\ref{eb5}) gives
\begin{eqnarray*}
\lefteqn{
\|\phi\|^2_{Y_1(T)} + \|\phi\|^2_{Z_1(T)}} \\
 &\leq& C(\|\nabla f\|^2_{L^2({\mathbb R}^3)} + \|g\|^2_{L^2({\mathbb
  R}^3)}) \nonumber \\
 & & + C\int_0^T \int_{{\mathbb R}^3} \left(
|\partial \phi||F| + \frac{|\phi||F|}{r^{1/2}\langle r\rangle^{1/2}} \right)dxdt
 \nonumber \\
 & & + C \left(\|h\|_{L^{\infty}(S_T)}
 + \|r^{1/2}\langle r \rangle^{1/2}\partial  h\|_{L^{\infty}(S_T)}
 \right)
 \nonumber \\
 & &\qquad \times
\Bigl(
\|r^{-5/4}\langle r\rangle^{-1/4} \phi\|^2_{L^2(S_T)}
+ \|r^{-1/4}\langle r\rangle^{-1/4} \partial \phi\|^2_{L^2(S_T)}\Bigr).
\nonumber
\end{eqnarray*}
Thus we immediately see the following lemma.


\begin{lemma}
 \label{ld7}
Let $\phi \in C^2(\overline{S_T}) \cap X_1(T)$ be a solution of the wave equation
\begin{equation}
 \label{edx9}
\partial_t^2 \phi - \Delta \phi + h(t,x)\Delta \phi = F
\quad \mbox{in}\,\, S_T,
\end{equation}
with initial data $(f,g) \in H^1({\mathbb R}^3)\times L^2({\mathbb R}^3)$.
Here $h \in C^1(\overline{S_T})$ is uniformly bounded
together with their first derivatives, and it also satisfies
\begin{eqnarray}
 \label{edx10}
 & & \|h\|_{L^{\infty}(S_T)} \leq \frac{1}{6}, \\
 \label{edx11}
 & & \|r^{1/2}\langle r \rangle^{1/2}\partial
 h\|_{L^{\infty}(S_T)} < \infty.
\end{eqnarray}
{\rm (i)} Suppose $r^{1/4}F \in L^2(S_T)$ and $0 \leq t \leq T$.
Then we have
\begin{eqnarray}
\label{edx12a}
\lefteqn{\|\partial \phi(t,\cdot)\|^2_{L^2({\mathbb R}^3)} +
\|\phi\|^2_{Y_1(T)} + \|\phi\|^2_{Z_1(T)}} \\
 &\leq& C_2(\|\nabla f\|^2_{L^2({\mathbb R}^3)} + \|g\|^2_{L^2({\mathbb
  R}^3)}) \nonumber \\
 & & + C_2(1+T)^{1/4}\|\phi\|_{Y_1(T)}\|r^{1/4}F\|_{L^2(S_T)} \nonumber \\
 & & + C_2(1+T)^{1/2}\|\phi\|^2_{Y_1(T)}
\begin{array}[t]{l}
 \bigl(\|h\|_{L^{\infty}(S_T)} \\
 \quad + \|r^{1/2}\langle r \rangle^{1/2}\partial  h\|_{L^{\infty}(S_T)}
 \bigr).
\end{array}
 \nonumber
\end{eqnarray}
{\rm (ii)} Suppose $r^{1/4}\langle r \rangle^{1/4}F \in L^2(S_T)$ and $0 \leq
 t \leq T$.
Then we have
\begin{eqnarray}
\label{edx12b}
\lefteqn{\|\partial \phi(t,\cdot)\|^2_{L^2({\mathbb R}^3)} +
\|\phi\|^2_{Y_1(T)} + \|\phi\|^2_{Z_1(T)}} \\
 &\leq& C_2(\|\nabla f\|^2_{L^2({\mathbb R}^3)} + \|g\|^2_{L^2({\mathbb
  R}^3)}) \nonumber \\
 & & + C_2 \bigl(\log (2+T)\bigr)^{1/2}
  \|\phi\|_{Z_1(T)}\|r^{1/4}\langle r \rangle^{1/4}F\|_{L^2(S_T)} \nonumber \\
 & & + C_2\log (2+T) \|\phi\|^2_{Z_1(T)}
\begin{array}[t]{l}
 \bigl(\|h\|_{L^{\infty}(S_T)} \\
 \quad + \|r^{1/2}\langle r \rangle^{1/2}\partial  h\|_{L^{\infty}(S_T)}
 \bigr).
\end{array}
 \nonumber
\end{eqnarray}
Here $C_2$ is a positive constant.
\end{lemma}

\subsection{Boundedness of $\{\phi_k\}$ in $X_2(T) \cap Y_2(T)$}


\begin{lemma}
 \label{ld1}
Let $(f,g) \in H^2_{{\rm rad}}({\mathbb R}^3) \times H^1_{{\rm rad}}({\mathbb R}^3)$
and set $\varepsilon = \|\nabla f\|_{H^1} + \|g\|_{H^1}$.
Then there exist positive constants $\varepsilon_1$, $A_1$
so that we can define the sequence $\{\phi_k\}$ {\rm (}see {\rm (\ref{ed14})},
{\rm (\ref{ed15}))} if $\varepsilon \leq \varepsilon_1$ and $T \leq \exp
 (A_1/\varepsilon)$.
In addition, there exists a constant $M_1$ such that
\begin{equation}
 \|\phi_k\|_{E_2(T)}
 + \|\phi_k\|_{Y_2(T)}
 + \|\phi_k\|_{Z_2(T)}
\leq M_1 \varepsilon
\quad \mbox{for}\,\, k = 0,1,2,\ldots
\nonumber
\end{equation}
Here $A_1$, $\varepsilon_1$, and $M_1$ depend on $a$, $b$ and $h$.
\end{lemma}
By Lemma \ref{ld2}, there exists a positive number
$c_0$ such that
\begin{equation}
 \label{ed33}
\|\phi\|_{E_2(T)} \leq c_0 \quad
\Longrightarrow \quad |h(\phi(t,x))| \leq \frac{1}{6}
\end{equation}
for $(t,x)\in \overline{S_T}$.
If $\|\phi_{k-1}\|_{E_2(T)} \leq c_0$,
then we can apply a standard existence theorem for
(\ref{ed14}), (\ref{ed15}) to define $\phi_k$.
We will see in the proof of Lemma \ref{ld1} that $\|\phi_k\|_{E_2(T)} \leq c_0$,
provided $\varepsilon \leq \varepsilon_1$ and $T \leq \exp (A_1/\varepsilon)$.
This ensures the well-definedness of the sequence $\{\phi_k\}$.


\begin{lemma}
 \label{ld6}
Let $\phi, {\tilde \phi} \in
C^{\infty}(\overline{S_T}) \cap X_2(T)$.
Assume that they satisfy
\begin{equation}
 \label{edx1}
 \partial_t^2 \phi - \Delta \phi + h({\tilde \phi})\Delta \phi
= F(\partial {\tilde \phi})
\quad \mbox{in}\,\,S_T
\end{equation}
and $\|{\tilde \phi}\|_{E_2(T)} \leq c_0$. Then we have
\begin{eqnarray}
 \label{edx2}
\lefteqn{\|\phi\|_{E_2(T)} + \|\phi\|_{Y_2(T)} + \|\phi\|_{Z_2(T)}} \\
&\leq& C_3\bigl(\|\partial \phi(0,\cdot)\|_{H^1({\mathbb R}^3)} +
 \|\partial {\tilde \phi}(0,\cdot)\|^2_{H^1 ({\mathbb R}^3)}\bigr) \nonumber \\
& & + C_3 \|{\tilde \phi}\|_{E_2(T)}
\bigl(\|\phi\|_{Z_2(T)} + \|{\tilde \phi}\|_{Z_2(T)}\bigr) \log
 (2 + T). \nonumber
\end{eqnarray}
Here, $C_3$ is a positive constant depending on $F$, $h$, $C_{\rm S}$
and $C_2$.
\end{lemma}

Proof of Lemma \ref{ld6}.
Differentiating (\ref{edx1}) with respect to $x_{\alpha}$ $(\alpha = 0,1,2,3)$, we see that
\begin{equation}
\label{edx13}
\partial^2_t (\partial_{\alpha} \phi) - \Delta (\partial_{\alpha} \phi) +
 h({\tilde \phi})\Delta (\partial_{\alpha} \phi)
= \partial_{\alpha}F(\partial {\tilde \phi})
- \partial_{\alpha} h({\tilde \phi})\Delta\phi.
\end{equation}
Since ${\tilde \phi}$ is radially symmetric by the definition of $X_2(T)$,
Lemma \ref{ld2} gives
\begin{eqnarray*}
 & &|F(\partial {\tilde \phi})|
\leq Cr^{-1/2}\langle r\rangle^{-1/2}\|{\tilde \phi}\|_{E_2(T)}|
\partial {\tilde \phi}|, \\
 & &|\partial_{\alpha}F(\partial {\tilde \phi})|
\leq Cr^{-1/2}\langle r\rangle^{-1/2}\|{\tilde \phi}\|_{E_2(T)}|
\partial^2 {\tilde \phi}|, \\
 & &|\partial h({\tilde \phi})|
 \leq Cr^{-1/2}\langle r \rangle^{-1/2} \|{\tilde \phi}\|_{E_2(T)}.
\end{eqnarray*}
Thus we get
\begin{eqnarray*}
& &\|r^{1/4}\langle r\rangle^{1/4}F(\partial {\tilde \phi})\|_{L^2_{t,x}} 
+ \|r^{1/4}\langle r\rangle^{1/4}\bigl(\partial_\alpha F(\partial {\tilde \phi})
 - \partial_\alpha h({\tilde \phi})\Delta \phi\bigr)\|_{L^2_{t,x}}
\\
& &\leq C\|{\tilde \phi}\|_{E_2(T)}
\bigl(\|{\tilde \phi}\|_{Z_2(T)}+\|\phi\|_{Z_2(T)}\bigr)
\bigl(\log (2+T)\bigr)^{1/2},
\end{eqnarray*}
where $L^p_{t,x}$ stands for $L^p(S_T)$.
Similarly, we have
\[
 \|h({\tilde \phi})\|_{L^{\infty}_{t,x}} + \|r^{1/2}\langle
 r \rangle^{1/2}\partial h({\tilde \phi})\|_{L^{\infty}_{t,x}}
 \leq C \|{\tilde \phi}\|_{E_2(T)}.
\]
Thanks to (\ref{ed33}), we also see that (\ref{edx10}) is satisfied for
$h({\tilde \phi})$.

Applying (\ref{edx12b}) to (\ref{edx1}) and (\ref{edx13}),
we obtain
\begin{eqnarray}
 \label{edx4}
 \lefteqn{\|\partial \phi(t,\cdot)\|^2_{L^2} +
 \|\partial^2 \phi(t,\cdot)\|^2_{L^2}
 + \|\phi\|^2_{Y_2(T)} + \|\phi\|^2_{Z_2(T)}}
\\
 &\leq & C\bigl(
\|\partial \phi(0,\cdot)\|^2_{L^2} + \|\partial^2 \phi(0,\cdot)\|^2_{L^2}
 \bigr)
 \nonumber \\
& &+C\|{\tilde \phi}\|_{E_2(T)}\|\phi\|_{Z_2(T)}
\bigl(\|{\tilde \phi}\|_{Z_2(T)}+\|\phi\|_{Z_2(T)}\bigr)
\log (2+T)
\nonumber \\
 &\leq & C\bigl(
\|\partial \phi(0,\cdot)\|^2_{L^2} + \|\partial^2 \phi(0,\cdot)\|^2_{L^2}
 \bigr) + \frac{1}{2}\|\phi\|^2_{Z_2(T)}
 \nonumber \\
& &
 + C\|{\tilde \phi}\|^2_{E_2(T)}
\bigl(\|{\tilde \phi}\|^2_{Z_2(T)}+\|\phi\|^2_{Z_2(T)}\bigr)
\bigl(\log (2+T)\bigr)^2
\nonumber
\end{eqnarray}
for $0 \leq t \leq T$.
In order to get (\ref{edx2}) from this estimate, we note that
\begin{eqnarray*}
 \|\partial_t^2 \phi(0,\cdot)\|_{L^2}
 &\leq& C\|\nabla^2 \phi(0,\cdot)\|_{L^2} + C\|\partial {\tilde \phi}(0,\cdot)\|^2_{L^4}
 \\
 &\leq& C\|\nabla^2 \phi(0,\cdot)\|_{L^2} + C\|\partial {\tilde \phi}(0,\cdot)\|_{H^1}^2
\end{eqnarray*}
by (\ref{edx1}) and Sobolev's inequality.
\hfill $\Box$

\bigskip

Proof of Lemma {\ref{ld1}}.
The proof proceeds by induction.
First, letting $\phi := \phi_0$, ${\tilde \phi} := 0$,
we apply Lemma \ref{ld6} to the equation
\begin{equation}
  \partial_t^2 \phi_0 - \Delta \phi_0 = 0
\nonumber
\end{equation}
and get
\begin{equation}
 \|\phi_0\|_{E_2(T)} + \|\phi_0\|_{Y_2(T)} + \|\phi_0\|_{Z_2(T)} \leq
C_3 (\|\nabla f_0\|_{H^1} + \|g_0\|_{H^1}).
\label{ed27}
\end{equation}
Since we easily see from the definition (\ref{ed11}) of $f_k$, $g_k$ that
\begin{equation}
\label{ed30}
\|\nabla f_k\|_{H^1} + \|g_k\|_{H^1}
  \leq \|\nabla f\|_{H^1} + \|g\|_{H^1}
 = \varepsilon
\end{equation}
for all $k = 0,1,2,\ldots$,
it follows from (\ref{ed27}) and (\ref{ed30}) with $k = 0$ that
\begin{equation}
 \|\phi_0\|_{E_2(T)} + \|\phi_0\|_{Y_2(T)} + \|\phi_0\|_{Z_2(T)}
\leq C_3 \varepsilon.
 \label{ed31}
\end{equation}
In order to define $\phi_1$ by (\ref{ed14})--(\ref{ed15}),
we need to verify that $\phi_0$ satisfies $\|\phi_0\|_{E_2(T)} \leq c_0$
(see (\ref{ed33})).
If we take $\varepsilon_1$ small so that the inequality $C_3\varepsilon_1 \leq c_0$ may hold, (\ref{ed31}) shows that it is true.
Now we are in a position to apply a standard existence, uniqueness, and
regularity theorem to (\ref{ed14})--(\ref{ed15}) with $k = 1$,
and we see that $\phi_1 \in C^{\infty}(\overline{S_T}) \cap X_2(T)$.
Moreover, by virtue of Lemma \ref{ld6}, $\phi_1$ enjoys the estimate
\begin{eqnarray}
 \label{ed32}
 \lefteqn{\|\phi_1\|_{E_2(T)} + \|\phi_1\|_{Y_2(T)} + \|\phi_1\|_{Z_2(T)}}
 \\
 &\leq & C_3 (\varepsilon + \varepsilon^2)
 \nonumber \\
 & & + C_3 \|\phi_0\|_{E_2(T)}
\bigl(\|\phi_1\|_{Z_2(T)}+\|\phi_0\|_{Z_2(T)}\bigr) \log (2+T),
 \nonumber
\end{eqnarray}
which immediately leads to
\begin{eqnarray}
 \label{ed42}
 \lefteqn{\|\phi_1\|_{E_2(T)} + \|\phi_1\|_{Y_2(T)} + \|\phi_1\|_{Z_2(T)}}
 \\
&\leq& \frac{4}{3}C_3 (\varepsilon + \varepsilon^2)
 + \frac{1}{3}\|\phi_0\|_{Z_2(T)},
\nonumber
\end{eqnarray}
provided that $T$ is chosen so that the inequality
$C_3\|\phi_0\|_{E_2(T)}\log (2+T) \leq 1/4$ may hold.
Thanks to the boundedness of $\|\phi_0\|_{E_2(T)}$ (see (\ref{ed31})),
we can choose $T$ in such a way that the last inequality holds.
Recalling that $\|\phi_0\|_{Z_2(T)}$ also has a similar bound,
we finally conclude that
\begin{equation}
\|\phi_1\|_{E_2(T)} + \|\phi_1\|_{Y_2(T)} + \|\phi_1\|_{Z_2(T)}
\leq M_1\varepsilon
\end{equation}
for $\varepsilon \leq \varepsilon_1$, where $M_1 = 4C_3$.

Before proceeding, we must take $\varepsilon_1$ still smaller so that
\begin{equation}
 M_1\varepsilon_1 \leq c_0
\quad \mbox{and}\quad
C_3M_1\varepsilon_1 \log 3 \leq \frac{1}{8}
\end{equation}
may hold.
Setting $A_1 := 1/(8C_3M_1)$,
we also take $T$ still smaller so that
$T \leq \exp (A_1/\varepsilon)$, thereby the inequality
\begin{eqnarray}
\label{ed43}
 \lefteqn{
C_3M_1\varepsilon \log (2+T)} \\
&\leq& C_3M_1\varepsilon \log (2+e^{A_1/\varepsilon}) \nonumber \\
&\leq& C_3M_1\varepsilon \log (3e^{A_1/\varepsilon}) \nonumber \\
&=& C_3M_1\varepsilon \log 3 + \frac{1}{8}
\leq \frac{1}{4} \nonumber
\end{eqnarray}
may hold for $\varepsilon \leq \varepsilon_1$.

Now we suppose that we have successively defined
$\phi_0,\,\phi_1,\,\ldots,\phi_{k-1} \in C^{\infty}(\overline{S_T}) \cap X_2(T)$
by (\ref{ed14})--(\ref{ed15}) satisfying
\begin{equation}
 \|\phi_\alpha\|_{E_2(T)} + \|\phi_\alpha\|_{Y_2(T)} + \|\phi_\alpha\|_{Z_2(T)}
 \leq M_1\varepsilon
\end{equation}
for $\alpha = 0,1,\ldots, k-1$.
A repetition of the above argument is all that is needed to show that,
according to (\ref{ed14})--(\ref{ed15}),
we can define $\phi_k$ in $C^{\infty}(\overline{S_T}) \cap X_2(T)$
satisfying the estimate
\begin{equation}
 \|\phi_k\|_{E_2(T)} + \|\phi_k\|_{Y_2(T)} + \|\phi_k\|_{Z_2(T)}
 \leq M_1\varepsilon.
\end{equation}
Lemma \ref{ld1} is thus proved by induction.
\hfill $\Box$

\subsection{Convergence in $X_1(T) \cap Y_1(T)$}


\begin{lemma}
 \label{ld9}
Let $\phi^{(i)}, {\tilde \phi}^{(i)} \in C^{\infty}(\overline{S_T}) \cap X_2(T)$ $(i = 1,2)$.
Assume that they satisty
\begin{equation}
 \partial_t^2 \phi^{(i)} - \Delta \phi^{(i)}
+ h({\tilde \phi}^{(i)})\Delta \phi^{(i)} =
 F(\partial {\tilde \phi^{(i)}})
\quad \mbox{in}\,\,S_T
\nonumber
\end{equation}
and $\|{\tilde \phi}^{(i)}\|_{E_2(T)} \leq c_0$ for $i = 1,2$. Then we have
\begin{eqnarray}
 \label{ed45}
\lefteqn{\|\phi^{(1)} - \phi^{(2)}\|_{E_1(T)}
+ \|\phi^{(1)} - \phi^{(2)}\|_{Y_1(T)}
+ \|\phi^{(1)} - \phi^{(2)}\|_{Z_1(T)}} \\
&\leq& C_4\|\partial \phi^{(1)}(0,\cdot)
- \partial \phi^{(2)}(0,\cdot)\|_{L^2}
\nonumber \\
& & + C_4 \bigl(\|{\tilde \phi}^{(1)}\|_{E_2(T)}
+ \|{\tilde \phi}^{(2)}\|_{E_2(T)}
+ \|\phi^{(2)}\|_{Y_2(T)}\bigr)(1+T)^{1/2}
 \nonumber \\
& & \quad \times \bigl(
\|{\tilde \phi}^{(1)} - {\tilde \phi}^{(2)}\|_{E_1(T)}
+ \|{\tilde \phi}^{(1)} - {\tilde \phi}^{(2)}\|_{Y_1(T)}
\nonumber \\
& & \qquad
+ \|\phi^{(1)} - \phi^{(2)}\|_{Y_1(T)}
\bigr).
\nonumber
\end{eqnarray}
Here, $C_4$ is a positive constant depending on $F$, $h$, $C_{\rm S}$, and $C_2$.
\end{lemma}

Proof. If we set $\phi^* = \phi^{(1)} - \phi^{(2)}$, it satisfies
\begin{eqnarray}
 \lefteqn{\partial_t^2 \phi^* - \Delta \phi^* +
 h({\tilde \phi}^{(1)})\Delta \phi^*}
\nonumber \\
&=&
F(\partial {\tilde \phi}^{(1)}) - F(\partial {\tilde \phi}^{(2)})
+\bigl(h({\tilde \phi}^{(2)}) - h({\tilde \phi}^{(1)})\bigr)\Delta
\phi^{(2)}. \nonumber
\end{eqnarray}
By (\ref{edx12a}), we have
\begin{eqnarray}
\lefteqn{\|\partial \phi^*(t,\cdot)\|^2_{L^2}
+ \|\phi^*\|^2_{Y_1(T)}+ \|\phi^*\|^2_{Z_1(T)}}
\nonumber\\
&\leq&\!\! C\|\partial \phi^*(0,\cdot)\|^2_{L^2} \nonumber\\
& &\!\! + C\|r^{1/4}\bigl(F(\partial {\tilde \phi}^{(1)}) - F(\partial
 {\tilde \phi}^{(2)})\bigr)\|_{L^2_{t,x}} \|\phi^*\|_{Y_1(T)}(1+T)^{1/4}
\nonumber \\
& &\!\! + C \|r^{1/4}\bigl(h({\tilde \phi}^{(2)})
- h({\tilde \phi}^{(1)})\bigr)\Delta
\phi^{(2)}\|_{L^2_{t,x}}\|\phi^*\|_{Y_1(T)}(1+T)^{1/4}
\nonumber \\
& &\!\! + C\bigl(
\|h({\tilde \phi}^{(1)})\|_{L^{\infty}_{t,x}}\!
+ \|r^{1/2}\langle r \rangle^{1/2}\partial
h({\tilde \phi}^{(1)})\|_{L^{\infty}_{t,x}}\bigr)
\|\phi^*\|^2_{Y_1(T)}(1+T)^{1/2}.
\nonumber
\end{eqnarray}
Here $L^p_{t,x}$ denotes $L^p(S_T)$.
By Lemma \ref{ld2}, we easily see that
\begin{eqnarray}
 \lefteqn{|F(\partial {\tilde \phi}^{(1)}) - F(\partial {\tilde \phi}^{(2)})|}
 \nonumber \\
 &\leq& Cr^{-1/2}\langle r\rangle^{-1/2}
(\|{\tilde \phi}^{(1)}\|_{E_2(T)}
+\|{\tilde \phi}^{(2)}\|_{E_2(T)})|\partial {\tilde \phi}^{(1)}
-\partial {\tilde \phi}^{(2)}|,
\nonumber \\
\lefteqn{|h({\tilde \phi}^{(2)})-h({\tilde \phi}^{(1)})|
 \leq Cr^{-1/2}\|{\tilde \phi}^{(2)}
- {\tilde \phi}^{(1)}\|_{E_1(T)},}
\nonumber \\
\lefteqn{\|h({\tilde \phi}^{(1)})\|_{L^{\infty}_{t,x}} +
\|r^{1/2}\langle r\rangle^{1/2}\partial
h({\tilde \phi}^{(1)})\|_{L^{\infty}_{t,x}}
 \leq C\|{\tilde \phi}^{(1)}\|_{E_2(T)}.}
\nonumber
\end{eqnarray}
Note that the assumption $\|{\tilde \phi}^{(i)}\|_{E_2(T)} \leq c_0$ is used here.
Comparing the second estimate with the first one,
we see that the factor $r^{-1/2}$ appearing in the second
estimate is weaker than the factor $r^{-1/2}\langle r\rangle^{-1/2}$ appearing
in the first one.
It is due to the second estimate that we have to use
(\ref{edx12a}) instead of (\ref{edx12b}).
Combining these estimates and the Schwarz inequality,
we obtain
\begin{eqnarray}
\lefteqn{\|\partial \phi^*(t,\cdot)\|^2_{L^2}
+ \|\phi^*\|^2_{Y_1(T)}+ \|\phi^*\|^2_{Z_1(T)}}
\nonumber \\
&\leq& C\|\partial \phi^*(0,\cdot)\|^2_{L^2}
\nonumber \\
& & + C\bigl(\|{\tilde \phi}^{(1)}\|_{E_2(T)}
+ \|{\tilde \phi}^{(2)}\|_{E_2(T)}\bigr)
\|\phi^*\|_{Y_1(T)}
\|{\tilde \phi}^*\|_{Y_1(T)}
(1+T)^{1/2} \nonumber \\
& & + C\|{\tilde \phi}^*\|_{E_1(T)}
\|\phi^*\|_{Y_1(T)}\|\phi^{(2)}\|_{Y_2(T)}
(1+T)^{1/2} \nonumber \\
& & + C\|{\tilde \phi}^{(1)}\|_{E_2(T)}
\|\phi^*\|^2_{Y_1(T)}(1+T)^{1/2} \nonumber \\
&\leq& C\|\partial \phi^* (0,\cdot)\|^2_{L^2}
+C\bigl(\|{\tilde \phi}^{(1)}\|_{E_2(T)} +\|{\tilde \phi}^{(2)}\|_{E_2(T)}
+ \|\phi^{(2)}\|_{Y_2(T)}\bigr)^2 \nonumber \\
& & \times (1+T)\bigl(\|{\tilde \phi}^*\|_{Y_1(T)} + \|{\tilde
 \phi}^*\|_{E_1(T)} + \|\phi^*\|_{Y_1(T)}\bigr)^2 + \frac{1}{2}\|\phi^*\|^2_{Y_1(T)},
\nonumber
\end{eqnarray}
where ${\tilde \phi}^* = {\tilde \phi}^{(1)} - {\tilde \phi}^{(2)}$.
This gives (\ref{ed45}).
\hfill $\Box$


\begin{lemma}
 \label{ld8}
Let $(f,g) \in H^2_{\rm rad}({\mathbb R}^3)\times H^1_{\rm rad}({\mathbb R}^3)$.
There exist positive constants $\varepsilon_2, A_2$ so that
the sequence $\{\phi_k\}$ defined in Lemma {\rm \ref{ld1}}
converges to a function $\phi \in X_1(T) \cap Y_1(T)$
if $\|\nabla f\|_{H^1} + \|g\|_{H^1} = \varepsilon \leq \varepsilon_2$
and $T \leq A_2 \varepsilon^{-2}$.
In addition, $\phi \in L^{\infty}([0,T]; H^2({\mathbb
 R}^3)) \cap C^{0,1}([0,T]; H^1({\mathbb
 R}^3)) \cap Y_2(T)$.
Here we denote by $C^{0,1}$ the space of Lipschitz continuous functions.
$A_2$ and $\varepsilon_2$ depend on $a$, $b$ and $h$.
\end{lemma}

Proof. Let $\varepsilon \leq \varepsilon_1$ and $T \leq \exp(A_1/\varepsilon)$.
Then $\{\phi_k\}$ is well-defined by Lemma \ref{ld1}.
Moreover, we have
\begin{equation}
 \|\phi_k\|_{E_2(T)}
 + \|\phi_k\|_{Y_2(T)}
 + \|\phi_k\|_{Z_2(T)}
\leq M_1 \varepsilon
\nonumber
\end{equation}
and $\|\phi_k\|_{E_2(T)}\leq c_0$
for $k = 0,1,2,\ldots$.
Therefore, if
\begin{equation}
 \label{ed52}
C_4\cdot 2M_1\varepsilon (1+T)^{1/2}
\leq \frac{1}{4},
\end{equation}
then it follows from Lemma \ref{ld9} that
\begin{eqnarray}
\lefteqn{\|\phi_k - \phi_{k-1}\|_{E_1(T)}
+ \|\phi_k - \phi_{k-1}\|_{Y_1(T)}
} \nonumber \\
&\leq& \frac{4}{3}C_4\|\partial \phi_k(0,\cdot)
- \partial \phi_{k-1}(0,\cdot)\|_{L^2}
 \nonumber \\
& & + \frac{1}{3}\bigl(\|\phi_{k-1} - \phi_{k-2}\|_{E_1(T)}
 + \|\phi_{k-1} - \phi_{k-2}\|_{Y_1(T)}
\bigr),
\nonumber
\end{eqnarray}
which, together with (\ref{ed12}), implies
\begin{equation}
 \sum_{k=1}^{\infty}\bigl(\|\phi_k - \phi_{k-1}\|_{E_1(T)}
+ \|\phi_k - \phi_{k-1}\|_{Y_1(T)}\bigr) < \infty.
\nonumber
\end{equation}
Hence, $\{\phi_k\}$ turns out to be a Cauchy sequence in $X_1(T)\cap Y_1(T)$
as long as $T \leq A_2 \varepsilon^{-2}$ and $\varepsilon \leq \varepsilon_2$,
where $A_2$ and $\varepsilon_2$ are chosen so that
\begin{equation}
 \label{ed54}
 \varepsilon_2 \leq \varepsilon_1,\quad
 2\sqrt{2}C_4 M_1\varepsilon_2 \leq \frac{1}{4},
\quad A_2 = \frac{1}{2\cdot 8^2C_4^2M_1^2}.
\end{equation}

Now it remains to prove that its limit $\phi$ belongs to
$L^{\infty}([0,T]; H^2({\mathbb R}^3)) \cap C^{0,1}([0,T]; H^1({\mathbb
R}^3)) \cap Y_2(T)$.
Since $\{\partial_{\alpha} \phi_k(t,\cdot)\}$ is bounded in
$H^1({\mathbb R}^3)$ and $\{\partial_{\alpha} \phi_{k}(t,\cdot)\}$ converges to
$\partial_{\alpha} \phi(t,\cdot)$ in $L^2({\mathbb R}^3)$, we see that
$\{\partial_{\alpha} \phi_{k}(t,\cdot)\}$ has a unique limit point
$\partial_{\alpha} \phi(t,\cdot)$ with respect to the weak topology of $H^1({\mathbb R}^3)$. Therefore,
\begin{equation}
 \label{ed71}
 \|\partial \phi(t,\cdot)\|_{H^1({\mathbb R}^3)} \leq \liminf_{k \to \infty}
 \|\partial \phi_{k}(t,\cdot)\|_{H^1({\mathbb R}^3)} \leq M_1 \varepsilon
\end{equation}
for $0 \leq t \leq T \leq A_2\varepsilon^{-2}$,
which shows that $\phi \in L^{\infty}([0,T]; H^2({\mathbb R}^3))
\cap C^{0,1}([0,T]; H^1({\mathbb R}^3))$.
By similar arguments, we conclude that
$\phi \in Y_2(T)$ and
\begin{equation}
 \|\phi\|_{Y_2(T)} + \|\phi\|_{Z_2(T)} \leq \liminf_{k \to \infty}
\bigl(\|\phi_k\|_{Y_2(T)} +  \|\phi_k\|_{Z_2(T)}\bigr) \leq M_1 \varepsilon.
\nonumber
\end{equation}
This completes the proof of Lemma \ref{ld8}.
\hfill $\Box$

\begin{remark}
\label{rd1}
Remark that (\ref{ed71}) holds for {\bf all} $t \in [0,T]$,
not to mention almost everywhere.
To be precise,
\begin{equation}
\|\partial \phi(t,\cdot)\|_{H^1({\mathbb R}^3)}
+ \|\phi\|_{Y_2(T)} + \|\phi\|_{Z_2(T)} \leq M_1\varepsilon
\quad \mbox{for all}\,\, t \in [0,T].
\nonumber
\end{equation}
Here $M_1$ is the constant appearing in Lemma \ref{ld1}.
\end{remark}


\section{Regularity}

\subsection{Preliminaries}
As a continuation of the previous section,
we will show that the limit $\phi$ of $\{\phi_k\}$ 
satisfies the equation (\ref{ea1}) in the strong sense.


\begin{lemma}
 \label{le3}
Let $(f,g) \in H^2_{\rm rad}({\mathbb R}^3)\times H^1_{\rm rad}({\mathbb R}^3)$
and $\varepsilon \leq \varepsilon_1$.
Assume that the sequence $\{\phi_k\}$ defined in Lemma {\rm \ref{ld1}}
converges to a function $\phi$ in $X_1(T) \cap Y_1(T)$
for some $T \leq \exp (A_1/\varepsilon)$.
Then $\phi \in X_2(T)$.
\end{lemma}

\noindent
As before, we set $\varepsilon = \|\nabla f\|_{H^1} + \|g\|_{H^1}$.
We have already proven in Lemma \ref{ld8} that
the assumptions in this lemma is satisfied
for $\varepsilon \leq \varepsilon_2$ and $T \leq A_2\varepsilon^{-2}$.
Note that the proof of Lemma \ref{ld8} implies $\|\partial \phi(t,\cdot)\|_{H^1({\mathbb R}^3)}
+ \|\phi\|_{Y_2(T)} + \|\phi\|_{Z_2(T)} \leq M_1\varepsilon$
for all $t \in [0,T]$ (see Remark \ref{rd1}).

We have to show that
$\partial_\alpha \partial_\beta \phi \in C^0([0,T]; L^2({\mathbb R}^3))$
for $\alpha, \beta = 0, 1, 2, 3$.
Let $t_0 \in [0,T]$, and let $\{t_n\}$ be an arbitrary sequence
such that
\[
 t_0 + t_n \in (0,T)\quad \mbox{for}\,\,n=1,2,3,\ldots
\quad \mbox{and} \quad
\lim_{n \to \infty} t_n = 0.
\]
We first aim at proving
\begin{equation}
\partial_a \partial_{\alpha}\phi (t_0 + t_n,\cdot) \to \partial_a
\partial_{\alpha} \phi(t_0,\cdot)\quad \mbox{in}\,\,L^2({\mathbb R}^3)
\quad \mbox{as}\,\,n \to \infty,
\nonumber
\end{equation}
for $a = 1,2,3$ and $\alpha = 0,1,2,3$.
In what follows, we set
\begin{equation}
 \psi(t,x) = \phi(t_0+t,x).
\nonumber
\end{equation}
By the assumption,
$\{\partial_{\alpha} \psi(t_n,\cdot)\}$
is a bounded sequence in $H^1({\mathbb R}^3)$.
Moreover,
$\{\partial_{\alpha}\psi(t_n,\cdot)\}$ converges to $\partial_{\alpha}
\psi(0,\cdot)$ in $L^2({\mathbb R}^3)$.
Hence, the sequence $\{\partial_a \partial_{\alpha} \psi(t_n,\cdot)\}$ converges
weakly to $\partial_a \partial_{\alpha} \psi(0,\cdot)$.
To prove the strong convergence,
we note the following fact: for a weakly convergent sequence $\{w_n\}$
in a Hilbert space $H$, the strong convergence is equivalent to
the estimate $\limsup_{n \to \infty} \|w_n\|_H \leq
\|w\|_H$.

In an abuse of notation, we set
\begin{equation}
\|[u,v]\|^2_{L^2({\mathbb R}^3)} = \int_{{\mathbb R}^3} \Bigl[
\frac{1}{2}|u(x)|^2 + \frac{1}{2}\bigl(1- h(\psi(0,x))\bigr)|v(x)|^2
\Bigr]\,dx
\nonumber
\end{equation}
for $u \in L^2({\mathbb R}^3)$, $v = (v_1,v_2,v_3) \in (L^2({\mathbb R}^3))^3$.
By (\ref{ed33}), this norm is equivalent to the usual norm of
$(L^2({\mathbb R}^3))^4$. So we aim at showing that
\begin{eqnarray}
 \label{ee4}
& & \limsup_{n \to \infty} \|[\partial_a \partial_t \psi(t_n,\cdot),
\partial_a \nabla \psi(t_n,\cdot)]\|_{L^2({\mathbb R}^3)} \\
& &\quad \leq \|[\partial_a \partial_t \psi(0,\cdot),
\partial_a \nabla \psi(0,\cdot)]\|_{L^2({\mathbb R}^3)}.\nonumber
\end{eqnarray}

Before showing (\ref{ee4}), let us prove two lemmas needed later.


\begin{lemma}
 \label{le1}
Let $\rho_k$ be the function defined by {\rm (\ref{ed7})}. Assume $0 \leq \alpha < 3$,
 $\beta \in {\mathbb R}$, and $0 \leq \gamma < 1$. Then we have
\begin{eqnarray}
 & &\int_{{\mathbb R}^3} \frac{\rho_k(y)}{|x - y|^\alpha
\langle x - y\rangle^\beta}dy \leq C|x|^{-\alpha}\langle x \rangle^{-\beta},
 \label{ee5a}\\
 & &\int_{{\mathbb R}^3} \frac{\rho_k(y)}{|x - y|^\alpha}dy
\leq C|x|^{-\alpha},
 \label{ee5}\\
& & \int_0^1 \frac{d\theta}{|x + \theta z|^{\gamma}\langle
x + \theta z\rangle^{\beta}} \leq C |x|^{-\gamma}\langle x \rangle^{-\beta}
\quad \mbox{for}\,\, |z| \leq 1,
 \label{ee6a} \\
& & \int_0^1 \frac{d\theta}{|x + \theta z|^{\gamma}}
\leq C |x|^{-\gamma}
\quad \mbox{for}\,\, |z| \leq 1.
 \label{ee6}
\end{eqnarray}
Here, in {\rm (\ref{ee5a})}--{\rm (\ref{ee5})}, $C$ is a constant independent of $k$.
\end{lemma}

Proof. We only prove (\ref{ee5}) and (\ref{ee6}),
since the others follow from similar arguments.
It is obvious that $(\ref{ee5})$ holds for $|x| \geq 3$,
because $2|x|/3 \leq |x - y| \leq 2|x|$ when $|y| \leq 1$.
So we suppose $|x|<3$ in what follows.
Since $|x|/2 \leq |x - y| \leq 3|x|/2$ for $|y| \leq |x|/2$,
we have
\[
 \int_{|y| \leq |x|/2} \frac{\rho_k(y)}{|x-y|^\alpha}\,dy
 \leq C|x|^{-\alpha}.
\]
If $|x|/2 < |y| < 1/k$, then we have $|x - y| < 3/k$ and $k < 2|x|^{-1}$.
Hence, we get
\begin{eqnarray*}
 \int_{|y| > |x|/2} \frac{\rho_k(y)}{|x-y|^\alpha}\,dy
 &\leq& \int_{|x - y| < 3/k} \frac{Ck^3}{|x - y|^{\alpha}}\,dy \\
 &\leq& Ck^{\alpha}
 \leq C|x|^{-\alpha}.
\end{eqnarray*}
Thus we have proven (\ref{ee5}).

We next prove (\ref{ee6}) only for $|x| < 3$, since
otherwise the inequality is obviously true.
If we set $x_z = x \cdot z / |z|$, then we see that
\begin{eqnarray}
 \label{ee7}
 \lefteqn{\int_0^1 |x + \theta z|^{-\gamma}d\theta} \\
&=& \int_0^1 (|x|^2 + 2\theta x_z |z| + \theta^2|z|^2)^{-\gamma/2}d\theta
\nonumber \\
&=& \frac{1}{|z|}\int_{x_z}^{x_z + |z|}(\lambda^2 +
 |x|^2 - x_z^2)^{-\gamma/2}d\lambda
\nonumber \\
&\leq& \frac{C}{|z|}\int_{x_z}^{x_z + |z|}\left(|\lambda| +
 \sqrt{|x|^2 - x_z^2}\right)^{-\gamma}d\lambda,
\nonumber
\end{eqnarray}
where $\lambda = x_z + \theta |z|$. Hence the case $|x_z| \leq |x|/2$ is
easily handled, because we can use $|\lambda | + \sqrt{|x|^2 - x_z^2} \geq
\sqrt{|x|^2 - (|x|/2)^2} = \sqrt{3}|x|/2$. So we assume $|x|/2 \leq
|x_z| \leq |x|$ in what follows. Let us set
\[
 \mathrm{I} = |z|^{-1}\int_{x_z}^{x_z + |z|}\left(|\lambda| +
 \sqrt{|x|^2 - x_z^2}\right)^{-\gamma}d\lambda.
\]
$\bullet$ case 1: $x_z \geq 0$.
\\
We simply estimate the integrand as $(\lambda + \sqrt{|x|^2 - x_z^2})^{-\gamma}
\leq \lambda^{-\gamma} \leq x_z^{-\gamma}$ and get
\begin{eqnarray*}
 \mathrm{I} &\leq& |z|^{-1}\int_{x_z}^{x_z + |z|}
  x_z^{-\gamma}\,d\lambda \\
 &\leq& |z|^{-1} (|x|/2)^{-\gamma}\int_{x_z}^{x_z + |z|} d\lambda
 = 2^{\gamma}|x|^{-\gamma}.
\end{eqnarray*}
$\bullet$ case 2: $x_z < 0 < x_z + |z|$.
Note that $|x_z| < |z|$ in this case.
\\
We can write the integral $\mathrm{I}$ as
$\mathrm{I} = \mathrm{I}_1 + \mathrm{I}_2$,
where
\begin{eqnarray*}
\mathrm{I}_1 &=& 2|z|^{-1}\int_{0}^{|x_z|}\left(\lambda +
 \sqrt{|x|^2 - x_z^2}\right)^{-\gamma}d\lambda, \\
\mathrm{I}_2 &=& |z|^{-1}\int_{|x_z|}^{x_z + |z|}\left(\lambda +
 \sqrt{|x|^2 - x_z^2}\right)^{-\gamma}d\lambda.
\end{eqnarray*}
Since $0 \leq \gamma < 1$ by assumption, we have
\begin{eqnarray*}
 \mathrm{I}_1 &\leq& 2|z|^{-1}\int_0^{|x_z|}\lambda^{-\gamma}\,d\lambda \\
 &=&2(1-\gamma)^{-1}|z|^{-1}|x_z|^{1-\gamma}
 \leq C|x|^{-\gamma}.
\end{eqnarray*}
We consider the estimate of $\mathrm{I}_2$ when $|x_z| \leq |z|/2$,
because $\mathrm{I}_2 < 0$ otherwise. We can treat this term as in case 1:
\[
 \mathrm{I}_2 \leq |z|^{-1}\int_{|x_z|}^{x_z + |z|} |x_z|^{-\gamma}\,d\lambda
 \leq C|x|^{-\gamma}.
\]
$\bullet$ case 3: $x_z + |z| \leq 0$. Note that $|z| \leq |x_z|$ in this case.
\\
By using the substitution $\lambda = -\mu$, we see
\[
 \mathrm{I} = |z|^{-1}\int_{|x_z|-|z|}^{|x_z|}\left(
\mu + \sqrt{|x|^2 - x_z^2}\right)^{-\gamma}\,d\mu.
\]
If $|z| \leq |x_z|/2$, then $|x_z| - |z| \geq |x_z|/2$. Thus, we get
\[
 \mathrm{I} \leq |z|^{-1}\int_{|x_z|-|z|}^{|x_z|} (|x_z|/2)^{-\gamma}\,d\mu
 \leq C|x|^{-\gamma}.
\]
If $|x_z|/2 \leq |z| \leq |x_z|$, then we estimate as
\[
 \mathrm{I} \leq |z|^{-1}\int_0^{|x_z|} \mu^{-\gamma}\,d\mu
 \leq (1-\gamma)^{-1}z|^{-1}|x_z|^{1-\gamma} \leq C|x|^{-\gamma}.
\]
Therefore, we have checked that (\ref{ee6}) holds for all the cases.
\hfill $\Box$


\begin{lemma}
\label{le2}
Let $\{\phi_k\}$ be the sequence defined in Lemma {\rm \ref{ld1}}.
Let $\phi$ be the function stated in Lemma {\rm \ref{le3}}.
Then for any $\eta \in {\mathbb R}^3$,
\begin{equation}
 \|r^{-1/4} \tau_\eta (\partial \phi_k -
 \partial \phi)\|_{L^2(S_T)} \to 0\quad
 \mbox{as}\,\, k \to \infty.
\nonumber
\end{equation}
Here, $\tau_\eta u(t,x) = u(t,x+\eta)$.
\end{lemma}

Proof.
Set $B(\eta) = \{x\, ; x \in {\mathbb R}^3,\ |x - \eta| < 1\}$ and $B(\eta)^c
= {\mathbb R}^3 \setminus B(\eta)$.
Then we can write
\begin{eqnarray}
\lefteqn{\|r^{-1/4}\tau_\eta(\partial \phi_k -
 \partial \phi)\|^2_{L^2(S_T)}}
\nonumber\\
&=& \|(\tau_{-\eta}r^{-1/4})(\partial \phi_k -
 \partial \phi)\|^2_{L^2(S_T)}
\nonumber \\
&=& \|(\tau_{-\eta}r^{-1/4})(\partial \phi_k -
 \partial \phi)\|^2_{L^2((0,T)\times B(\eta))}
\nonumber \\
& &+ \|(\tau_{-\eta}r^{-1/4})(\partial \phi_k -
 \partial \phi)\|^2_{L^2((0,T)\times B(\eta)^c)}.
\nonumber 
\end{eqnarray}
It is easy to see that $\lim_{k \to \infty} \|(\tau_{-\eta}r^{-1/4})(\partial \phi_k -
 \partial \phi)\|^2_{L^2((0,T)\times B(\eta)^c)} = 0$, because
$\|(\tau_{-\eta}r^{-1/4})(\partial \phi_k -
 \partial \phi)\|^2_{L^2((0,T)\times B(\eta)^c)} \leq T \|\phi_k - \phi\|^2_{E_1(T)}$.
To handle the other term, we use Lemma \ref{ld2} and get
\begin{eqnarray*}
 |\partial \phi_k(t,x) - \partial \phi(t,x)|
 &\leq& C|x|^{-1/2}(\|\phi_k\|_{E_2(T)} + \|\phi\|_{E_2(T)}) \\
 &\leq& C|x|^{-1/2}.
\end{eqnarray*}
This gives
\begin{eqnarray*}
\lefteqn{\|(\tau_{-\eta}r^{-1/4})(\partial \phi_k -
 \partial \phi)\|^2_{L^2((0,T)\times B(\eta))}
} \\
 &\leq& C \int_0^T \int_{|x-\eta|<1} |x-\eta|^{-1/2}|x|^{-1/2}|\partial \phi_k(t,x) -
 \partial \phi(t,x)|dxdt,
\end{eqnarray*}
and hence $\lim_{k \to \infty} \|(\tau_{-\eta}r^{-1/4})(\partial \phi_k -
 \partial \phi)\|^2_{L^2((0,T)\times B(\eta))} = 0$ by the Schwarz inequality.
\hfill $\Box$

\subsection{Proof of Lemma \ref{le3}}

We set
\begin{equation}
 \psi_k(t,x) = \phi_{k}(t_0+t,x)\quad
\mbox{for}\,\,t_0+t\in (0,T),\ k=0,1,2,\ldots
\nonumber
\end{equation}
Then $\psi_k$ satisfies
\begin{eqnarray}
 \label{ee10}
\lefteqn{\partial_t^2 (\Delta_{\eta} \psi_k) - \bigl(1 - h(\psi_{k-1})\bigr)
\Delta (\Delta_{\eta}\psi_k)} \\
& & =
\Delta_{\eta} F(\partial \psi_{k-1})
-\Delta_{\eta} h(\psi_{k-1})\Delta
\tau_{\eta} \psi_k
\nonumber
\end{eqnarray}
for $\eta \in {\mathbb R}^3,\,\, k=0,1,2,\ldots$,
where
\begin{equation}
 \tau_{\eta} u(t,x) = u(t,x+ \eta),\quad
\Delta_{\eta} = \tau_{\eta} - 1.
\nonumber
\end{equation}
Integrating (\ref{ee10})$\times \partial_t \Delta_{\eta} \psi_k$ as usual, we have
\begin{eqnarray}
\lefteqn{\int_{{\mathbb R}^3} \Bigl[\frac{1}{2}(\partial_t \Delta_{\eta}
\psi_k)^2 + \frac{1}{2}\bigl(1-h(\psi_{k-1})\bigr)|\nabla \Delta_{\eta} \psi_k|^2
\Bigr](t,x)\,dx}
\nonumber \\
&=&\int_{{\mathbb R}^3} \Bigl[\frac{1}{2}(\partial_t \Delta_{\eta}
\psi_k)^2 + \frac{1}{2}\bigl(1-h(\psi_{k-1})\bigr)|\nabla \Delta_{\eta} \psi_k|^2
\Bigr](0,x)\,dx
\nonumber \\
& &+\int_0^t\!\! \int_{{\mathbb R}^3} \Bigl[
\partial_t \Delta_{\eta} \psi_k
 \bigl(\Delta_{\eta} F(\partial \psi_{k-1})
- \Delta_{\eta} h(\psi_{k-1}) \Delta \tau_{\eta} \psi_k\bigr) \nonumber \\
& &\qquad - \frac{1}{2}\partial_t h(\psi_{k-1})|\nabla \Delta_{\eta}
 \psi_k|^2
+\partial_a h(\psi_{k-1})\partial_t\Delta_{\eta}\psi_k \partial^a
 \Delta_{\eta} \psi_k \Bigr]\,dxds. \nonumber
\end{eqnarray}
Note that
\begin{eqnarray*}
 \lefteqn{|h(\psi_k(t,x)) - h(\psi(t,x))|} \\
 &\leq& C\|\psi_k(t,\cdot) - \psi(t,\cdot)\|^{1/4}_{L^2}
\|\psi_k(t,\cdot) -
 \psi(t,\cdot)\|^{3/4}_{W^{1,6}} \\
 &\leq& C\|\psi_k(t,\cdot) - \psi(t,\cdot)\|^{1/4}_{L^2}
\|\nabla \psi_k(t,\cdot) -
 \nabla \psi(t,\cdot)\|^{3/4}_{H^1}
\end{eqnarray*}
by the Gagliardo-Nirenberg inequality.
Hence, passing to the limit,
we immediately see that
\begin{eqnarray}
 \label{ee13}
\lefteqn{\int_{{\mathbb R}^3} \Bigl[\frac{1}{2}(\partial_t \Delta_{\eta}
\psi)^2 + \frac{1}{2}\bigl(1-h(\psi)\bigr)|\nabla \Delta_{\eta} \psi|^2
\Bigr](t,x)\,dx} \\
&\leq& \|[\partial_t \Delta_{\eta} \psi(0,\cdot),
 \nabla \Delta_{\eta} \psi(0,\cdot)]\|^2_{L^2({\mathbb R}^3)}
+ C\limsup_{k \to \infty}({\rm I}_{k}
+ {\rm II}_{k} + {\rm III}_{k}),
\nonumber
\end{eqnarray}
where we set
\begin{eqnarray}
& &{\rm I}_k=\int_{I_t}\int_{{\mathbb R}^3} |\partial_t
\Delta_\eta \psi_k||\partial \Delta_\eta \psi_{k-1}|
\bigl(|\tau_\eta \partial \psi_{k-1}| + |\partial \psi_{k-1}|\bigr)\,dxds, \\
& &{\rm II}_k=\int_{I_t}\int_{{\mathbb R}^3} |\partial_t
\Delta_\eta \psi_k||\Delta_\eta \psi_{k-1}||\Delta \tau_\eta\psi_k|\,dxds,
\nonumber\\
& &{\rm III}_k=\int_{I_t}\int_{{\mathbb R}^3} |\partial
 \psi_{k-1}| |\partial \Delta_\eta \psi_k|^2\,dxds
\nonumber
\end{eqnarray}
and
\begin{equation}
 I_t = \left\{
\begin{array}{ll}
(0,t) & \mbox{if}\,\, t>0, \\
(t,0) & \mbox{if}\,\, t<0.
\end{array}
\right.
\end{equation}

Let us consider ${\rm II}_{k}$ first. By Lemma \ref{ld2} and
the Schwarz inequality,
we have
\begin{eqnarray}
{\rm II}_{k}\!\!&=&\!\!\int_{I_t}\int_{{\mathbb R}^3} |\partial_t
\Delta_{-\eta} \psi_{k}||\Delta_{-\eta}\psi_{k-1}||\Delta \psi_{k}|\,dxds
\nonumber \\
&\leq&\!\! C\|\Delta_{-\eta}\psi_{k-1}\|_{L^{\infty}(I_t;
 H^1)}\int_{I_t}\int_{{\mathbb R}^3} \frac{|\partial_t
\Delta_{-\eta} \psi_{k}||\Delta \psi_{k}|}{r^{1/2}}\,dxds
\nonumber \\
&\leq&\!\! C\|\Delta_{-\eta}\psi_{k-1}\|_{L^{\infty}(I_t; H^1)}
\|r^{-1/4}\Delta_{-\eta}\partial_t\psi_{k}\|_{L^2(I_t\times {\mathbb R}^3)}
\nonumber \\
& & \times
\|r^{-1/4}\Delta \psi_{k}\|_{L^2(I_t\times {\mathbb R}^3)}.
\nonumber
\end{eqnarray}
Hence, noting
Lemma \ref{le2} and Lemma \ref{ld1}, we get
\begin{equation}
\label{ee17}
 \limsup_{k \to \infty} {\rm II}_{k} \leq
 C\|\Delta_{-\eta}\psi\|_{L^{\infty}(I_t; H^1)}
\|r^{-1/4}\Delta_{-\eta}\partial_t \psi\|_{L^2(I_t\times {\mathbb R}^3)}.
\end{equation}
As for ${\rm III}_{k}$, we estimate it as
\begin{equation}
{\rm III}_{k}
 \leq C\|\partial \psi_{k-1}\|_{L^{\infty}(I_t; H^1)}\|r^{-1/4} \partial
 \Delta_{\eta} \psi_{k}\|^2_{L^2(I_t \times {\mathbb R}^3)}
\nonumber
\end{equation}
by Lemma \ref{ld2}.
Noting Lemma \ref{le2} and Lemma \ref{ld1} again,
we see
\begin{equation}
 \label{ee18}
\limsup_{k \to \infty} {\rm III}_{k}
 \leq C\|r^{-1/4}\Delta_{\eta} \partial \psi\|^2_{L^2(I_t \times {\mathbb R}^3)}.
\end{equation}
Similarly,
\begin{eqnarray}
 \label{ee18a}
\lefteqn{\limsup_{k \to \infty} {\rm I}_{k}} \\
 &\leq& C\|r^{-1/4}\Delta_{\eta}\partial \psi\|^2_{L^2(I_t \times {\mathbb R}^3)}
+C\|r^{-1/4}\Delta_{-\eta}\partial \psi\|^2_{L^2(I_t \times {\mathbb R}^3)}.
\nonumber
\end{eqnarray}
As a consequence of (\ref{ee13}) and (\ref{ee17})--(\ref{ee18a}),
we obtain
\begin{eqnarray}
 \label{ee19}
\lefteqn{\int_{{\mathbb R}^3} \Bigl[\frac{1}{2}(\partial_t \Delta_{\eta}
\psi)^2 + \frac{1}{2}\bigl(1-h(\psi)\bigr)|\nabla \Delta_{\eta} \psi|^2
\Bigr](t,x)\,dx} \\
&\leq& \|[\partial_t \Delta_{\eta} \psi(0,\cdot), \nabla \Delta_{\eta}
 \psi(0,\cdot)]\|^2_{L^2({\mathbb R}^3)}
\nonumber \\
& & + C\|\Delta_{-\eta}\psi\|_{L^{\infty}(I_t;
 H^1)}\|r^{-1/4}\Delta_{-\eta}\partial_t \psi\|_{L^2(I_t \times {\mathbb R}^3)}
\nonumber \\
& & + C\|r^{-1/4}\Delta_{\eta}\partial \psi\|^2_{L^2(I_t \times {\mathbb R}^3)}
+ C\|r^{-1/4}\Delta_{-\eta}\partial \psi\|^2_{L^2(I_t \times {\mathbb R}^3)}.
\nonumber
\end{eqnarray}

Now we set $\eta = \delta e_j$, where $e_j$ is the unit coordinate vector
in the $x_j$ direction. Then it is well known that
\begin{eqnarray}
 \label{ee20}
 & & \int_{{\mathbb R}^3} \frac{1}{2}(\partial_t \Delta_{\delta e_j} \psi(0,x))^2dx
 \leq \delta^2 \int_{{\mathbb R}^3} \frac{1}{2}(\partial_t \partial_j
 \psi(0,x))^2dx, \\
 \label{ee21}
 & & \|\Delta_{-\delta e_j}\psi\|_{L^{\infty}(I_t; H^1)} \leq
 |\delta|\|\partial_j \psi\|_{L^{\infty}(I_t; H^1)}.
\end{eqnarray}
In order to estimate $\int_{{\mathbb R}^3}\frac{1}{2}\bigl(1-h(\psi(0,x))\bigr)|\nabla
\Delta_{\delta e_j}\psi(0,x)|^2dx$, we note that
\[
\|\partial_a \Delta_{\delta e_j} \psi(0,\cdot) - \delta \partial_a
 \partial_j \psi(0,\cdot)\|_{L^2} = o(\delta),
\]
where $o(\delta)/\delta \to 0$ as $\delta \to 0$.
Therefore, in view of
\begin{eqnarray*}
 \lefteqn{
|\nabla \Delta_{\delta e_j} \psi(0,x)|^2 - \delta^2 |\nabla \partial_j \psi(0,x)|^2
} \\
&=& \bigl(\partial_a \Delta_{\delta e_j} \psi(0,x) - \delta \partial_a
 \partial_j \psi(0,x)\bigr)\partial^a \Delta_{\delta e_j} \psi(0,x) \\
& & + \delta\partial_a\partial_j \psi(0,x)\bigl(\partial^a
 \Delta_{\delta e_j} \psi(0,x) - \delta \partial^a
 \partial_j \psi(0,x)
\bigr),
\end{eqnarray*}
we see that
\begin{eqnarray}
 \label{ee22}
\lefteqn{
\int_{{\mathbb R}^3}\frac{1}{2}\bigl(1-h(\psi(0,x))\bigr)|\nabla
\Delta_{\delta e_j}\psi(0,x)|^2dx
}\\
&=& \delta^2 \int_{{\mathbb R}^3}\frac{1}{2}\bigl(1-h(\psi(0,x))\bigr)|\nabla
\partial_j \psi(0,x)|^2dx + o(\delta^2), \nonumber
\end{eqnarray}
where $o(\delta^2)/\delta^2 \to 0$ as $\delta \to 0$.
Hence from (\ref{ee19})--(\ref{ee22}), we have
\begin{eqnarray}
 \label{ee23}
\lefteqn{\int_{{\mathbb R}^3} \Bigl[\frac{1}{2}(\partial_t
\Delta_{\delta e_j} \psi)^2 + \frac{1}{2}\bigl(1-h(\psi)\bigr)|\nabla
\Delta_{\delta e_j} \psi|^2\Bigr](t,x)\,dx} \\
&\leq& \delta^2 \|[\partial_t \partial_j \psi(0,\cdot), \nabla \partial_j
 \psi(0,\cdot)]\|^2_{L^2({\mathbb R}^3)}
\nonumber \\
& & + C|\delta|\|\partial_j\psi\|_{L^{\infty}(I_t;
 H^1)}\|r^{-1/4}\Delta_{-\delta e_j}\partial_t \psi\|_{L^2(I_t \times {\mathbb R}^3)}
\nonumber \\
& & + C\|r^{-1/4}\Delta_{\delta e_j}\nabla \psi\|^2_{L^2(I_t \times
 {\mathbb R}^3)} \nonumber \\
& & + C\|r^{-1/4}\Delta_{-\delta e_j}\nabla \psi\|^2_{L^2(I_t \times {\mathbb R}^3)}
 + o(\delta^2).
\nonumber
\end{eqnarray}

It remains to estimate $\|r^{-1/4}\Delta_{\delta e_j}\partial_{\alpha}
\psi\|_{L^2(I_t \times {\mathbb R}^3)}$ $(\alpha = 0,1,2,3)$.
We first notice that
\begin{equation}
 \label{ee24}
 \|r^{-1/4}\Delta_{\delta e_j}\partial_{\alpha}\psi\|_{L^2(I_t \times
 {\mathbb R}^3)}
= \lim_{l \to \infty}  \|r^{-1/4}\rho_l *\Delta_{\delta
e_j}\partial_{\alpha}\psi\|_{L^2(I_t \times {\mathbb R}^3)}.
\end{equation}
In fact, we can apply a similar argument as Lemma \ref{le2} to the
estimate of $\|r^{-1/4}\bigl(\rho_l *\Delta_{\delta
e_j}\partial_{\alpha}\psi -
\Delta_{\delta e_j}\partial_{\alpha}\psi\bigr)\|^2_{L^2(I_t \times
 {\mathbb R}^3)}$
by dividing the domain of integration as ${\mathbb R}^3 = \{r < 1\} \cup
\{r \geq 1\}$.
Since $\rho_l * \Delta_{\delta e_j}\partial_{\alpha}\psi(s,x)$ is smooth
with respect to $x$, we have
\[
 |\rho_l * \Delta_{\delta e_j}\partial_{\alpha}\psi(s,x)|^2
\leq \delta^2 \int_0^1 |\rho_l * \partial_j \partial_{\alpha}\psi(s,x +
\theta \delta e_j)|^2d\theta,
\]
hence we get
\begin{eqnarray*}
\lefteqn{\|r^{-1/4}\rho_l * \Delta_{\delta
 e_j}\partial_{\alpha}\psi(s,\cdot)\|^2_{L^2({\mathbb R}^3)}} \\
&\leq& \delta^2 \int_0^1\int_{{\mathbb R}^3}r^{-1/2}|\rho_l * \partial_j \partial_{\alpha}\psi(s,x +
\theta \delta e_j)|^2dx d\theta.
\nonumber
\end{eqnarray*}
Furthermore, the Schwarz inequality gives
\[
 |\rho_l * \partial_j \partial_{\alpha} \psi(s,x + \theta\delta e_j)|^2
 \leq \int_{{\mathbb R}^3}\rho_l(x + \theta\delta e_j - y)|\partial_j
 \partial_{\alpha} \psi(s,y)|^2\,dy.
\]
Thus we obtain
\begin{eqnarray}
 \label{ee25}
 \lefteqn{\|r^{-1/4}\rho_l * \Delta_{\delta
 e_j}\partial_{\alpha}\psi(s,\cdot)\|^2_{L^2({\mathbb R}^3)}} \\
&\leq& \delta^2 \int_{{\mathbb R}^3} a_l(y,\delta)
|\partial_j  \partial_{\alpha} \psi(s,y)|^2\,dy,
\nonumber
\end{eqnarray}
where
\begin{equation}
 a_l(y,\delta) = \int_0^1 \int_{{\mathbb R}^3}
|x|^{-1/2}\rho_l(x + \theta\delta e_j - y)\,dxd\theta
\quad (|\delta| \leq 1).
\nonumber
\end{equation}
By (\ref{ee5}) and (\ref{ee6}), we have
\begin{equation}
 \label{ee27}
 a_l(y,\delta) \leq C\int_0^1 |\theta \delta e_j - y|^{-1/2}d\theta
 \leq C|y|^{-1/2}.
\end{equation}

Consequently, as a result of (\ref{ee24})--(\ref{ee27}), we obtain
\begin{equation}
 \label{ee28}
 \|r^{-1/4}\Delta_{\delta e_j}\partial_{\alpha}\psi\|_{L^2(I_t \times
 {\mathbb R}^3)} \leq C|\delta| \|r^{-1/4}\partial_j
 \partial_{\alpha}\psi\|_{L^2(I_t \times {\mathbb R}^3)}.
\end{equation}
Therefore, by (\ref{ee23}) and (\ref{ee28}), we arrive at
\begin{eqnarray}
\lefteqn{\int_{{\mathbb R}^3} \Bigl[\frac{1}{2}(\partial_t
\Delta_{\delta e_j} \psi)^2 + \frac{1}{2}\bigl(1-h(\psi)\bigr)|\nabla
\Delta_{\delta e_j} \psi|^2\Bigr](t,x)\,dx}
\nonumber \\
&\leq& \delta^2 \|[\partial_t \partial_j \psi(0,\cdot), \nabla \partial_j
 \psi(0,\cdot)]\|^2_{L^2({\mathbb R}^3)}
\nonumber \\
& & + C\delta^2\|\partial_j\psi\|_{L^{\infty}(I_t;
 H^1)}\|r^{-1/4}\partial_j \partial_t \psi\|_{L^2(I_t \times {\mathbb R}^3)}
\nonumber \\
& & + C\delta^2 \|r^{-1/4}\partial_j \partial \psi\|^2_{L^2(I_t \times {\mathbb R}^3)}
 + o(\delta^2),
\nonumber
\end{eqnarray}
which further yields
\begin{eqnarray}
 \label{ee30}
\lefteqn{\int_{{\mathbb R}^3} \Bigl[\frac{1}{2}(\partial_t
\partial_j \psi)^2 + \frac{1}{2}\bigl(1-h(\psi)\bigr)|\nabla
\partial_j \psi|^2\Bigr](t,x)\,dx} \\
&\leq& \|[\partial_t \partial_j \psi(0,\cdot), \nabla \partial_j
 \psi(0,\cdot)]\|^2_{L^2({\mathbb R}^3)}
\nonumber \\
& & + C\|\partial_j\psi\|_{L^{\infty}(I_t;
 H^1)}\|r^{-1/4}\partial_j \partial_t \psi\|_{L^2(I_t \times {\mathbb R}^3)}
+ C\|r^{-1/4}\partial_j \partial \psi\|^2_{L^2(I_t \times {\mathbb R}^3)}
\nonumber
\end{eqnarray}
for $t_0 + t \in (0,T)$.

Finally, noting
\begin{eqnarray}
 \lefteqn{|h\bigl(\psi(0,x)\bigr) - h\bigl(\psi(t,x)\bigr)|}
 \nonumber \\
 &\leq& C\|\psi(0,\cdot) - \psi(t,\cdot)\|^{1/4}_{L^2}
 \|\nabla \psi(0,\cdot) - \nabla \psi(t,\cdot)\|^{3/4}_{H^1},
\nonumber
\end{eqnarray}
(\ref{ee30}) leads to (\ref{ee4}).
Thus we have proven
\begin{equation}
 \label{ee31}
 \partial_a \partial_{\alpha} \phi \in C^0([0,T];
 L^2({\mathbb R}^3))\quad
\mbox{for}\,\, a = 1,2,3,\ \alpha = 0,1,2,3.
\end{equation}
To show that $\partial^2_t \phi \in
C^0([0,T]; L^2({\mathbb R}^3))$, we use the equation (\ref{ed14}).
Let $v \in C_0^{\infty}((0,T)\times {\mathbb R}^3)$.
Since the sequence $\{\phi_k\}$ converges to $\phi$ in $X_1(T) \cap Y_1(T)$,
it is easy to verify that the limit of the equation
\begin{eqnarray*}
& &-\int_0^T \int_{{\mathbb R}^3} \bigl(
\partial_t \phi_k \partial_t v - \nabla \phi_k \cdot \nabla v + \nabla
\phi_k \cdot \nabla (h(\phi_{k-1})v)\bigr)dxdt \\
& &\quad = \int_0^T \int_{{\mathbb R}^3} F(\partial \phi_{k-1}) v dxdt
\end{eqnarray*}
as $k \to \infty$ gives
\begin{eqnarray*}
& &-\int_0^T \int_{{\mathbb R}^3} \bigl(
\partial_t \phi \partial_t v - \nabla \phi \cdot \nabla v + \nabla
\phi \cdot \nabla (h(\phi)v)\bigr)dxdt \\
& &\quad = \int_0^T \int_{{\mathbb R}^3} F(\partial \phi) v dxdt.
\end{eqnarray*}
Therefore, we obtain $\partial^2_t \phi = \Delta \phi - h(\phi)\Delta
\phi + F(\partial \phi)$ in $(0,T) \times {\mathbb R}^3$.
Since $H^2({\mathbb R}^3) \subset L^\infty({\mathbb R}^3)$ and
$H^1({\mathbb R}^3) \subset L^4({\mathbb R}^3)$, this implies
$\partial^2_t \phi \in C^0([0,T]; L^2({\mathbb R}^3))$. \hfill $\Box$


\section{Almost Global Existence}

In this section we complete the proof of Theorem \ref{ta1}.
So far, we have proven that there exists a strong solution for
the initial value problem (\ref{ea1}), (\ref{ea2}) for $T \leq
A_2\varepsilon^{-2}$,
where $\varepsilon = \|\nabla f\|_{H^1} + \|g\|_{H^1}$ and $\varepsilon
\leq \varepsilon_2$.
It remains to show that the solution actually extends
to the time $t \leq \exp (A_1\varepsilon^{-1})$.

Let us define a sequence $\{\psi^{T_1}_k\}$ by (\ref{ed14}), (\ref{ed15}),
where we now set $f_k$, $g_k$ as
\begin{equation}
 f_k = \rho_{2^k}*\phi(T_1,\cdot),\quad
 g_k = \rho_{2^k}*\partial_t\phi(T_1,\cdot)
\nonumber
\end{equation}
for $T_1 = A_2\varepsilon_2^{-2}$.
Note that $\|\nabla \phi(T_1,\cdot)\|_{H^1}
+ \|\partial_t \phi(T_1,\cdot)\|_{H^1} \leq M_1\varepsilon$ (see Remark \ref{rd1}).
Suppose that $\varepsilon$ is small so that $M_1\varepsilon \leq \varepsilon_2$.
If we choose $T$ smaller so that $T \leq A_2(M_1\varepsilon)^{-2}$,
then we see by Lemma \ref{ld1} and Lemma \ref{ld8} that
$\|\psi_k^{T_1}\|_{E_2(T)} + \|\psi_k^{T_1}\|_{Y_2(T)}
+ \|\psi_k^{T_1}\|_{Z_2(T)} \leq M_1^2\varepsilon$,
and $\{\psi^{T_1}_k\}$ converges to a function $\psi^{T_1}$ in
$X_1(T) \cap Y_1(T)$.
In addition, $\psi^{T_1} \in X_2(T) \cap Y_2(T)$ and
$\|\psi^{T_1}\|_{E_2(T)} + \|\psi^{T_1}\|_{Y_2(T)}
+ \|\psi^{T_1}\|_{Z_2(T)} \leq M_1^2\varepsilon$.
We consider another sequence $\{\phi_k^{T_1}\}$ on $[0,T]\times {\mathbb R}^3$,
defined by $\phi^{T_1}_k(t,x) = \phi_k(t+T_1,x)$.
Then it follows from Lemma \ref{ld9} that
\begin{eqnarray}
 \lefteqn{\|\psi_k^{T_1} - \phi_k^{T_1}\|_{E_1(T)}
+ \|\psi_k^{T_1} - \phi_k^{T_1}\|_{Y_1(T)}}
\nonumber \\
&\leq& C_4\|\partial \psi^{T_1}_k(0,\cdot)
- \partial \phi^{T_1}_k(0,\cdot)\|_{L^2} \nonumber \\
& & + C_4(\|\psi^{T_1}_{k-1}\|_{E_2(T)}
+ \|\phi^{T_1}_{k-1}\|_{E_2(T)}
+ \|\phi^{T_1}_{k}\|_{Y_2(T)})(1+T)^{1/2} \nonumber\\
& & \times \bigl(
\|\psi_{k-1}^{T_1} - \phi_{k-1}^{T_1}\|_{E_1(T)}
+ \|\psi_{k-1}^{T_1} - \phi_{k-1}^{T_1}\|_{Y_1(T)}
+ \|\psi_k^{T_1} - \phi_k^{T_1}\|_{Y_1(T)}
\bigr).
\nonumber
\end{eqnarray}
We take $A_3$ and $\varepsilon_0$ so that
\begin{eqnarray}
 & &M_1 \varepsilon_0 \leq \varepsilon_2,\quad
C_4(2M_1\varepsilon_0 + M_1^2\varepsilon_0)
\cdot 2^{1/2}\leq \frac{1}{4}, \nonumber\\
& & A_3 \leq A_2M_1^{-2},\quad
A_3 \leq \frac{1}{2\cdot 4^2C_4^2M_1^2(2+M_1)^2}.
\nonumber
\end{eqnarray}
Then for $T_2 = A_3\varepsilon^{-2}_0$ and $\varepsilon \leq \varepsilon_0$, we have
\begin{equation}
 C_4 (2M_1 \varepsilon + M_1^2 \varepsilon) (1+T_2)^{1/2} \leq \frac{1}{4},
\nonumber
\end{equation}
and thus
\begin{eqnarray}
 \lefteqn{\|\psi^{T_1}_k - \phi^{T_1}_k\|_{E_1(T_2)}
+ \|\psi^{T_1}_k - \phi^{T_1}_k\|_{Y_1(T_2)}}
\nonumber \\
 &\leq& \frac{4}{3}C_4\|\partial \psi^{T_1}_k(0,\cdot)
- \partial \phi^{T_1}_k(0,\cdot)\|_{L^2}
 \nonumber \\
 & & + \frac{1}{3}\bigl(\|\psi^{T_1}_{k-1} - \phi^{T_1}_{k-1}\|_{E_1(T_2)}
+ \|\psi^{T_1}_{k-1} - \phi^{T_1}_{k-1}\|_{Y_1(T_2)}\bigr).
\nonumber
\end{eqnarray}
Hence we have $\|\psi^{T_1}_{k} - \phi^{T_1}_{k}\|_{E_1(T_2)}
+ \|\psi^{T_1}_{k} - \phi^{T_1}_{k}\|_{Y_1(T_2)} \to 0$ as $k \to
\infty$, which implies $\phi^{T_1}_k \to \psi^{T_1}$
in $X_1(T_2) \cap Y_1(T_2)$.
Therefore, we have shown that
there exists an extension of $\phi \in X_1(T_1) \cap
Y_1(T_1)$ such that
$\phi_k \to \phi$ in $X_1(T_1 + T_2)
\cap Y_1(T_1 + T_2)$, and it turns out that
$\|\phi\|_{E_2(T_1 + T_2)}
+ \|\phi\|_{Y_2(T_1 + T_2)}
+ \|\phi\|_{Z_2(T_1 + T_2)} \leq M_1\varepsilon$.
Repeating this argument, we can extend the solution as
long as $T \leq \exp
(A_1 \varepsilon^{-1})$.
This completes the proof of Theorem \ref{ta1}.


\section{Continuous dependence on the data}

\subsection{Preliminaries}
We have proven the almost global existence of a
strong solution to the initial value problem (\ref{ea1}), (\ref{ea2}).
In this section, we further discuss the question whether the solution
depends continuously on the data (Theorem \ref{ta2}).
In order to prove Theorem \ref{ta2},
we should weaken the regularity assumptions of Lemma \ref{ld7}
so as to apply the estimates (\ref{edx12a})--(\ref{edx12b})
to low-regularity solutions.


\begin{lemma}
 \label{lg1}
Let $\phi \in C^2(\overline{S_T}) \cap X_1(T)$
be a solution of {\rm (\ref{edx9})} satisfying
$\|\phi\|_{C^2(\overline{S_T})} < \infty$
with initial data $(f,g) \in H^1({\mathbb R}^3) \times L^2({\mathbb R}^3)$.
Assume {\rm (\ref{edx10})}, {\rm (\ref{edx11})} and $h \in X_1(T)$.
Then the assertions {\rm (i)}, {\rm (ii)} in Lemma {\rm \ref{ld7}} remain valid.
\end{lemma}

\bigskip

Proof. We set
\begin{equation}
\label{eg1}
 h_k(t,\cdot) = \rho_k * h(t,\cdot)
\quad \mbox{for}\,\, k = 1,2,3,\ldots,
\end{equation}
which is the convolution of $\rho_k$ and $h(t,\cdot)$ as functions
on ${\mathbb R}^3$.
Recalling the definition of $\rho_k$,
we see that $h_k \in C^1(\overline{S_T})$ and $\|h_k\|_{C^1(\overline{S_T})} < \infty$.
Since $\phi$ satisfies the equation
\begin{equation}
 \label{eg2}
\partial_t^2 \phi - \Delta \phi + h_k(t,x)\Delta \phi
= F + (h_k - h) \Delta \phi,
\end{equation}
we have,
applying the inequality preceding Lemma \ref{ld7},
\begin{eqnarray}
\lefteqn{
\|\partial \phi(t,\cdot)\|^2_{L^2({\mathbb R}^3)}
+\|\phi\|^2_{Y_1(T)} + \|\phi\|^2_{Z_1(T)}
}\nonumber \\
&\leq& C\bigl(\|\nabla \phi(0,\cdot)\|^2_{L^2} +
 \|\partial_t \phi(0,\cdot)\|^2_{L^2}\bigr) \nonumber \\
& & + C\int_0^T \int_{{\mathbb R}^3}\Bigl(
|\partial \phi||F| + \frac{|\phi||F|}{r^{1/2}\langle r\rangle^{1/2}} \Bigr)dxdt \nonumber \\
& & + C\int_0^T \int_{{\mathbb R}^3}\Bigl(
|\partial \phi||(h_k - h) \Delta \phi|
+ \frac{|\phi||(h_k - h) \Delta \phi|}{r^{1/2}\langle r\rangle^{1/2}} \Bigr)dxdt \nonumber \\
& & + C\bigl(\|h_k\|_{L^{\infty}_{t,x}} +
 \|r^{1/2}\langle r\rangle^{1/2}\partial h_k\|_{L^{\infty}_{t,x}}\bigr)
 \nonumber \\
& &\qquad \times \bigl(
\|r^{-5/4}\langle r\rangle^{-1/4}\phi\|^2_{L^2_{t,x}} +
\|r^{-1/4}\langle r \rangle^{-1/4}\partial \phi\|^2_{L^2_{t,x}}
\bigr) \nonumber
\end{eqnarray}
for $0 \leq t \leq T$.
Here $L^p_{t,x}$ stands for $L^p(S_T)$.
Using Lemma \ref{le1}, we have
\begin{eqnarray}
 \label{eg4}
 |\partial h_k(t,x)| &\leq& \int_{{\mathbb R}^3} \rho_k(y)|\partial
 h(t,x-y)|\,dy \\
&\leq& \|r^{1/2}\langle r\rangle^{1/2}\partial h\|_{L^{\infty}_{t,x}}
 \int_{{\mathbb R}^3}\frac{\rho_k(y)}{|x-y|^{1/2}\langle
 x-y\rangle^{1/2}}dy \nonumber \\
&\leq& C\|r^{1/2}\langle r\rangle^{1/2}\partial
 h\|_{L^{\infty}_{t,x}}|x|^{-1/2}\langle x\rangle^{-1/2}, \nonumber
\end{eqnarray}
and hence
\[
\|h_k\|_{L^{\infty}_{t,x}} +
 \|r^{1/2}\langle r\rangle^{1/2}\partial h_k\|_{
L^{\infty}_{t,x}}
\leq C\bigl(\|h\|_{L^{\infty}_{t,x}} +
\|r^{1/2}\langle r\rangle^{1/2}\partial h\|_{L^{
\infty}_{t,x}} \bigr).
\]
Therefore, noting
\begin{eqnarray*}
\lefteqn{
\int_0^T \int_{{\mathbb R}^3}\Bigl(
|\partial \phi||(h_k - h) \Delta \phi|
+ \frac{|\phi||(h_k - h) \Delta \phi|}{r^{1/2}\langle r\rangle^{1/2}} \Bigr)dxdt
} \\
&\leq&C\|\Delta \phi\|_{L^{\infty}_{t,x}}
\bigl(\|\partial \phi\|_{L^2_{t,x}}
+ \|\phi/r\|_{L^2_{t,x}}
\bigr) \|h_k - h\|_{L^2_{t,x}}
\end{eqnarray*}
and passing to the limit, we obtain the
desired estimate. \hfill $\Box$


\begin{lemma}
 \label{lg2}
Let $\phi \in X_1(T) \cap Y_1(T)$ be a weak solution of {\rm (\ref{edx9})}
with initial data $(f,g) \in H^1({\mathbb R}^3) \times L^2({\mathbb R}^3)$.
Assume {\rm (\ref{edx10})}, {\rm (\ref{edx11})} and $h \in X_1(T)$.
Then the assertions {\rm (i)}, {\rm (ii)} in Lemma {\rm \ref{ld7}} remain valid.
\end{lemma}
Here the constant $C_2$ appearing in (\ref{edx12a})--(\ref{edx12b})
should be replaced by a larger one,
which is written as $C_2$ again.
We say that $\phi \in X_1(T) \cap Y_1(T)$ is a weak solution of
(\ref{edx9}) if
\begin{equation}
 \label{egx1}
 -\int_0^T\!\! \int_{{\mathbb R}^3} (\partial_t \phi \partial_t \psi -
  \partial^j \phi \partial_j \psi + \partial^j \phi \partial_j(h\psi)) dxdt
 = \int_0^T\!\! \int_{{\mathbb R}^3} F \psi dxdt
\end{equation}
holds for all $\psi \in C_0^{\infty}(S_T)$.

\bigskip

Proof.
Let ${\bar \phi}$ be the zero-extension of $\phi$, i.e.
\begin{equation}
 \label{eg5}
{\bar \phi}(t,x) = \left\{
\begin{array}{ll}
 \phi(t,x)& \mbox{for}\,\,0 \leq t \leq T, \\
 0 & \mbox{otherwise}.
\end{array}
\right.
\end{equation}
We also extend $h$, $F$ to ${\bar h}$, ${\bar F}$ similarly.
We set
\begin{equation}
 {\bar \rho}_k(t,x) = \rho^{(1)}_k\left(t + \frac{1}{k}\right)\rho_k(x)
\quad \mbox{for}\,\, (t,x)\in {\mathbb R}\times {\mathbb R}^3,
\nonumber
\end{equation}
where $\rho^{(1)}_k(t)$ is a 1-dimensional counterpart of (\ref{ed7}).
We define
\begin{equation}
 {\bar \phi}_k(t,x) = {\bar \rho}_k \bar{*} {\bar \phi}(t,x),
\nonumber
\end{equation}
where $u {\bar *} v$ denotes the convolution of functions $u, v$ on ${\mathbb R} \times {\mathbb R}^3$.
Here we take $k$ large enough, so that $1/k < T/4$.

We first show that
\begin{eqnarray}
 \label{eg11}
\lefteqn{
\partial_t^2 {\bar \phi}_k - \Delta {\bar \phi}_k
+ h(t,x)\Delta {\bar \phi}_k} \\
&=& {\bar \rho}_k {\bar *} {\bar F} + {\bar h} \Delta {\bar \phi}_k
 - \partial^j {\bar \rho}_k {\bar *} ({\bar h} \overline{\partial_j \phi}) +
 {\bar \rho}_k {\bar *} (\overline{\partial^j h}\, \overline{\partial_j \phi})
\nonumber
\end{eqnarray}
for $0 < t < T - 2/k$. Here, $\overline{\partial_j \phi}$ and
$\overline{\partial^j h}$ denote the zero-extensions of $\partial_j \phi$
and $\partial^j h$, respectively. To verify this, we
compute $\langle \partial_t^2 {\bar \phi}_k - \Delta {\bar \phi}_k, \psi
\rangle$ for $\psi \in C_0^{\infty}({\mathbb R} \times {\mathbb R}^3)$,
assuming $\mathrm{supp}\, \psi \subset (0, T-2/k) \times {\mathbb R}^3$.
Note that if we define a product $\tilde{*}$ by
\[
 u \tilde{*} v = \int_{{\mathbb R}} \int_{{\mathbb R}^3} u(s -
 t, y - x) v(s, y) dy ds,
\]
then $\mathrm{supp}\, {\bar \rho}_k \tilde{*} \psi \subset (0,
T) \times {\mathbb R}^3$. Therefore, we easily see that
\[
 \langle \partial_{\alpha}^2 {\bar \phi}_k, \psi \rangle
 = - \int_0^T\!\! \int_{{\mathbb R}^3} \partial_{\alpha} \phi \,
 \partial_{\alpha} ({\bar \rho}_k \tilde{*} \psi) dxdt.
\]
Hence, by (\ref{egx1}), we get
\begin{eqnarray*}
 \lefteqn{\langle \partial_t^2 {\bar \phi}_k - \Delta {\bar \phi}_k,
  \psi \rangle} \\
 &=& \int_0^T\!\! \int_{{\mathbb R}^3} F\cdot
 {\bar \rho}_k \tilde{*} \psi dxdt
 + \int_0^T\!\! \int_{{\mathbb R}^3} \partial^j \phi\,
 \partial_j(h\cdot {\bar \rho}_k \tilde{*} \psi) dxdt \\
 &=& \langle {\bar \rho}_k {\bar *} {\bar F} + {\bar \rho}_k {\bar *}
  (\overline{\partial^j \phi}\,\overline{\partial_j h}) - \partial_j
  {\bar \rho}_k {\bar *} ({\bar h} \overline{\partial^j \phi}), \psi \rangle.
\end{eqnarray*}
This gives (\ref{eg11}).

We next show that
\begin{eqnarray}
 \label{eg12}
\lefteqn{
\|r^{1/4}\langle r\rangle^{1/4}[h\Delta {\bar \phi}_k
- \partial^j {\bar \rho}_k {\bar *} ({\bar h} \overline{\partial_j \phi})]\|_{L^2
((0,T - 2/k)\times {\mathbb R}^3)}} \\
\lefteqn{
+ \|r^{1/4}\langle r\rangle^{1/4}
{\bar \rho}_k {\bar *} (\overline{\partial^j h}\, \overline{\partial_j \phi})\|_{L^2
((0,T - 2/k)\times {\mathbb R}^3)}} \nonumber \\
&\leq&\!\!\!\! C\|r^{1/2}\langle r\rangle^{1/2}\partial h\|_{L^{\infty}
(S_T)}
\|r^{-1/4}\langle r\rangle^{-1/4}\nabla \phi\|_{L^2
(S_T)} \nonumber
\end{eqnarray}
for $k = 1,2,3,\ldots$.
For this purpose, set $u_j = \partial_j \phi$ and
\begin{eqnarray}
 & &{\rm I}_k= \int_0^T \int_{{\mathbb R}^3}
 \partial^j {\bar \rho}_k(t-s,x-y)
 (h(t,x) - h(s,y)) u_j(s,y)\,dyds,
\nonumber \\
 & & {\rm II}_k= \int_0^T \int_{{\mathbb R}^3}
 {\bar \rho}_k(t-s,x-y)\partial^j h(s,y)u_j(s,y)\,dyds.
\nonumber
\end{eqnarray}
In what follows, we will prove
\begin{eqnarray*}
\lefteqn{
\|r^{1/4}\langle r\rangle^{1/4} {\rm I}_k\|_{L^2 (S_{T_k})}
+ \|r^{1/4}\langle r\rangle^{1/4} {\rm II}_k\|_{L^2 (S_{T_k})}} \\
&\leq&\!\!\!\! C\|r^{1/2}\langle r\rangle^{1/2}\partial h\|_{L^{\infty}
(S_T)}
\|r^{-1/4}\langle r\rangle^{-1/4}\nabla \phi\|_{L^2
(S_T)},
\end{eqnarray*}
where we set $S_{T_k} = (0,T-2/k)\times {\mathbb R}^3$.

We first consider the estimate on ${\rm II}_k$.
Since
\begin{eqnarray*}
 |{\rm II}_k| &=& \left|\int_0^T \rho^{(1)}_k\left(t - s + \frac{1}{k}\right)
\rho_k * (\partial^j h(s,\cdot)u_j(s,\cdot))ds\right| \\
 &\leq& \left(
\int_0^T \rho^{(1)}_k\left(t - s + \frac{1}{k}\right)
\left|\rho_k * (\partial^j h(s,\cdot)u_j(s,\cdot))\right|^2ds
\right)^{1/2}
\end{eqnarray*}
by the Schwarz inequality, we have
\begin{eqnarray}
\lefteqn{
\|r^{1/4}\langle r \rangle^{1/4}{\rm II}_k\|^2_{L^2(S_{T_k})}}
\nonumber \\
&\leq & \int_0^T \int_{{\mathbb R}^3}|x|^{1/2}\langle x\rangle^{1/2}
\left|\rho_k * (\partial^j h(s,\cdot)u_j(s,\cdot))\right|^2dxds
\nonumber \\
&\leq & \|r^{1/2}\langle r\rangle^{1/2}\partial
 h\|^2_{L^{\infty}(S_T)}
\int_0^T\int_{{\mathbb R}^3}|x|^{1/2}\langle x\rangle^{1/2}
\nonumber \\
& &\qquad \qquad\times
\left(\int_{{\mathbb R}^3}\frac{\rho_k(x-y)|\nabla \phi(s,y)|}{|y|^{1/2}\langle
 y\rangle^{1/2}}dy\right)^2dxds. \nonumber
\end{eqnarray}
To estimate the last integral, let us set
\begin{eqnarray}
& &{\rm II}_k^{(1)} = \int_0^T\int_{{\mathbb R}^3}|x|^{1/2}\langle x\rangle^{1/2}
\nonumber \\
& &\qquad \qquad \quad
 \times\left(\int_{|y|<1/k}\frac{\rho_k(x-y)|\nabla \phi(s,y)|}{|y|^{1/2}\langle
 y\rangle^{1/2}}dy\right)^2dxds, \nonumber\\
& &{\rm II}_k^{(2)} = \int_0^T\int_{{\mathbb R}^3}|x|^{1/2}\langle x\rangle^{1/2}
\nonumber \\
& &\qquad \qquad \quad
 \times\left(\int_{|y|>1/k}\frac{\rho_k(x-y)|\nabla \phi(s,y)|}{|y|^{1/2}\langle
 y\rangle^{1/2}}dy\right)^2dxds. \nonumber
\end{eqnarray}
Since $|x| < 2/k$ if $|x-y|<1/k$ and $|y|<1/k$, it follows
\begin{eqnarray*}
 \lefteqn{
\int_{{\mathbb R}^3}|x|^{1/2}\langle x\rangle^{1/2}
\int_{|y|<1/k}\frac{\rho_k(x-y)^2}{|y|^{1/2}}dydx} \\
&\leq& C\int_{{\mathbb R}^3} k^{-1/2}
\int_{|y|<1/k}\frac{k^3\rho_k(x-y)}{|y|^{1/2}}dydx \\
&=& C \int_{|y|<1/k} \frac{k^{5/2}}{|y|^{1/2}} dy
 \leq C.
\end{eqnarray*}
Therefore, by the Schwarz inequality, we obtain
\begin{eqnarray}
{\rm II}_k^{(1)} &\leq& \int_0^T\int_{{\mathbb R}^3}|x|^{1/2}\langle x\rangle^{1/2}
\left(\int_{|y|<1/k}\frac{\rho_k(x-y)^2}{|y|^{1/2}}dy\right)
\nonumber \\
& &\quad \times
\|r^{-1/4}\langle r\rangle^{-1/4}\nabla \phi(s,\cdot)\|^2_{L^2}dxds \nonumber \\
&\leq& C\|r^{-1/4}\langle r\rangle^{-1/4}\nabla
 \phi\|^2_{L^2(S_T)}. \nonumber
\end{eqnarray}
Meanwhile, since $|x| \leq 2|y|$ if $|x-y|<1/k$ and $|y|>1/k$, we have
\begin{eqnarray}
{\rm II}_k^{(2)}&\leq& C\int_0^T \int_{{\mathbb R}^3}
\left(\int_{{\mathbb R}^3}\frac{\rho_k(x-y)|\nabla \phi(s,y)|}{|y|^{1/4}\langle
 y\rangle^{1/4}}dy\right)^2dxds
\nonumber \\
&\leq& C\int_0^T \int_{{\mathbb R}^3} \int_{{\mathbb R}^3}
 \frac{\rho_k(x-y)|\nabla \phi(s,y)|^2}{|y|^{1/2}\langle
 y\rangle^{1/2}} dy dx ds \nonumber \\
&\leq& C \|r^{-1/4}\langle r\rangle^{-1/4}\nabla
 \phi\|^2_{L^2(S_T)}. \nonumber
\end{eqnarray}
Thus we conclude that
\begin{eqnarray}
\lefteqn{\|r^{1/4}\langle r \rangle^{1/4}{\rm II}_k\|_{L^2(S_{T_k})}}
\nonumber \\
&\leq&\!\!\!\! C\|r^{1/2}\langle r\rangle^{1/2}\partial
 h\|_{L^{\infty}(S_T)}
 \|r^{-1/4}\langle r\rangle^{-1/4}\nabla
 \phi\|_{L^2(S_T)}. \nonumber
\end{eqnarray}

We next turn ourself to the estimate of ${\rm I_k}$.
Suppose $|t-s|\leq 1/k$, $|x-y|\leq 1/k$ and set $z = y - x$.
Using Lemma \ref{le1}, we get
\begin{eqnarray*}
\lefteqn{|h(t,x) - h(s,y)|}\\
&\leq& |(t-s,x-y)|\int_0^1 |\partial h(t + \theta (s-t), x +
 \theta (y-x))|\,d\theta \\
&\leq& Ck^{-1}\|r^{1/2}\langle r\rangle^{1/2}\partial
 h\|_{L^{\infty}(S_T)}
\int_0^1 \frac{d\theta}{|x + \theta z|^{1/2}\langle x +
\theta z\rangle^{1/2}} \\
&\leq& Ck^{-1}\|r^{1/2}\langle r\rangle^{1/2}\partial
 h\|_{L^{\infty}(S_T)}|x|^{-1/2}\langle x\rangle^{-1/2}
\end{eqnarray*}
for $0 \leq t \leq T - 2/k$.
Strictly speaking, here we should first approximate ${\bar h}$
as (\ref{eg1}), and then use (\ref{eg4}).
This justification is easy, so we omit it.

It follows from the above estimate that
\begin{eqnarray*}
|{\rm I}_k| &\leq& C\|r^{1/2}\langle r\rangle^{1/2}\partial
h\|_{L^{\infty}(S_T)}|x|^{-1/2}\langle
x\rangle^{-1/2} \\
& & \times \int_0^T\!\! \int_{{\mathbb R}^3}
k^{-1} |\partial^j {\bar \rho}_k(t-s,x-y)||u_j(s,y)|\,dyds. \nonumber
\end{eqnarray*}
Since $\int_0^T\int_{{\mathbb R}^3}k^{-1}|\partial^j {\bar
\rho}_k(t-s,x-y)|\,dyds \leq C$, we proceed as
\begin{eqnarray*}
|{\rm I}_k|&\leq& C\|r^{1/2}\langle r\rangle^{1/2}\partial
h\|_{L^{\infty}(S_T)}|x|^{-1/2}\langle
x\rangle^{-1/2} \nonumber\\
& & \times \left(\int_0^T\!\! \int_{{\mathbb R}^3}
 k^{-1}|\partial^j {\bar \rho}_k(t-s,x-y)||u_j(s,y)|^2\,dyds\right)^{1/2} \nonumber
\end{eqnarray*}
for $0 \leq t \leq T - 2/k$. This gives
\begin{eqnarray}
 \lefteqn{\|r^{1/4}\langle r\rangle^{1/4}{\rm I}_k\|^2_{L^2(S_{T_k})}
\leq C\|r^{1/2}\langle r\rangle^{1/2}\partial
h\|^2_{L^{\infty}(S_T)}
}\nonumber \\
& & \qquad \qquad \qquad \times \int_0^T\!\!\int_{{\mathbb R}^3}\!\!
 \int_{{\mathbb R}^3}
 \frac{k^{-1}|\partial^j \rho_k(x-y)||u_j(s,y)|^2}{|x|^{1/2}\langle
 x\rangle^{1/2}}\,dydxds. \nonumber
\end{eqnarray}
Since we are able to show
\[
 \int_{{\mathbb R}^3} k^{-1}\frac{|\partial^j
 \rho_k(x-y)|}{|x|^{1/2}\langle x \rangle^{1/2}}dx \leq
 C|y|^{-1/2}\langle y \rangle^{-1/2}
\]
in the same way as in the proof of Lemma \ref{le1}, we arrive at
\begin{eqnarray}
\lefteqn{\|r^{1/4}\langle r \rangle^{1/4}{\rm I}_k\|_{L^2(S_{T_k})}}
\nonumber \\
&\leq&\!\!\!\! C\|r^{1/2}\langle r\rangle^{1/2}\partial
 h\|_{L^{\infty}(S_T)}
 \|r^{-1/4}\langle r\rangle^{-1/4}\nabla
 \phi\|_{L^2(S_T)}. \nonumber
\end{eqnarray}
This completes the proof of (\ref{eg12}).

Now we are in a position to apply Lemma \ref{lg1} to (\ref{eg11}).
We can see that $\|{\bar \phi}_k\|_{C^2(S_T)}
< \infty$ by the Schwarz inequality.
Let us prove below that (\ref{edx12b}) holds in the present setting if
$r^{1/4}\langle r \rangle^{1/4}F \in L^2(S_T)$.
The inequality (\ref{edx12a}) can be proven similarly. We set
\[
 F_k = {\bar \rho}_k {\bar *} {\bar F} + {\bar h} \Delta {\bar \phi}_k
 - \partial^j {\bar \rho}_k {\bar *} ({\bar h} \overline{\partial_j \phi}) +
 {\bar \rho}_k {\bar *} (\overline{\partial^j h}\, \overline{\partial_j \phi})
\]
(the right-hand side of (\ref{eg11})) and estimate $\|r^{1/4}\langle
r\rangle^{1/4}F_k\|_{L^2(S_{T_k})}$.
By using (\ref{eg12}), we get
\begin{eqnarray}
 \label{egx2}
 \lefteqn{\|r^{1/4}\langle r\rangle^{1/4}F_k\|_{L^2(S_{T_k})}
 \leq \|r^{1/4}\langle r\rangle^{1/4}{\bar \rho}_k {\bar *} {\bar F}\|_{L^2(S_{T_k})}}
  \\
& &\qquad + C\|r^{1/2}\langle r\rangle^{1/2}\partial
 h\|_{L^{\infty}(S_T)}
 \|r^{-1/4}\langle r\rangle^{-1/4}\nabla
 \phi\|_{L^2(S_T)}.
\nonumber
\end{eqnarray}
Furthermore, we see that
\[
 \|r^{1/4}\langle r\rangle^{1/4}{\bar \rho}_k {\bar *} {\bar F}\|_{L^2(S_{T_k})}
 \leq C \|r^{1/4}\langle r\rangle^{1/4}F\|_{L^2(S_T)}.
\]
In fact, if we set $G(t,x) = r^{1/4}\langle r\rangle^{1/4}F(t,x)$, we have
\begin{eqnarray*}
 |{\bar \rho}_k {\bar *} {\bar F}(t,x)|^2
 \!\!&=&\!\! \left|
\int_0^T\int_{{\mathbb R}^3}\frac{{\bar
\rho}_k(t-s,x-y)G(s,y)}{|y|^{1/4}\langle y \rangle^{1/4}}dyds
\right|^2\\
 &\leq&\!\! Cr^{-1/2}\langle r\rangle^{-1/2}\!\!\int_0^T\!\!\int_{{\mathbb R}^3}\!
 {\bar \rho}_k(t-s,x-y)|G(s,y)|^2dyds
\end{eqnarray*}
by the Schwarz inequality and Lemma \ref{le1},
which immediately yields
\begin{eqnarray*}
 \lefteqn{\|r^{1/4}\langle r\rangle^{1/4}{\bar \rho}_k {\bar *} {\bar
  F}\|^2_{L^2(S_{T_k})}} \\
 &\leq& C \int_0^{T-2/k}\!\!\!\int_{{\mathbb R}^3}\int_0^{T}\!\int_{{\mathbb R}^3}
 {\bar \rho}_k(t-s,x-y)|G(s,y)|^2dydsdxdt \\
 &\leq& C\|G\|^2_{L^2(S_T)} = \|r^{1/4}\langle r\rangle^{1/4}F\|^2_{L^2(S_T)},
\end{eqnarray*}
as desired.
Thus (\ref{egx2}) gives
\begin{eqnarray}
 \lefteqn{\|r^{1/4}\langle r\rangle^{1/4}F_k\|_{L^2(S_{T_k})}
 \leq C\|r^{1/4}\langle r\rangle^{1/4} F\|_{L^2(S_T)}}
  \\
& &\qquad + C\|r^{1/2}\langle r\rangle^{1/2}\partial
 h\|_{L^{\infty}(S_T)}
 \|r^{-1/4}\langle r\rangle^{-1/4}\nabla
 \phi\|_{L^2(S_T)}.
\nonumber
\end{eqnarray}
Therefore, by Lemma \ref{lg1}, we have
\begin{eqnarray}
\label{egx3}
\lefteqn{
\|\partial {\bar \phi}_k(t,\cdot)\|^2_{L^2({\mathbb R}^3)}
+\|{\bar \phi}_k\|^2_{Y_1(T-2/k)} + \|{\bar \phi}_k\|^2_{Z_1(T-2/k)}
} \\
 &\leq&
C\bigl(\|\nabla {\bar \phi}_k(0,\cdot)\|^2_{L^2} +
 \|\partial_t {\bar \phi}_k(0,\cdot)\|^2_{L^2}\bigr)
\nonumber \\
 & & + C \bigl(\log (2+T)\bigr)^{1/2}
  \|{\bar \phi}_k\|_{Z_1(T-2/k)}
\nonumber \\
 & & \quad \times
\begin{array}[t]{l}
\bigl(
\|r^{1/4}\langle r\rangle^{1/4} F\|_{L^2(S_T)}\\
 + \|r^{1/2}\langle r\rangle^{1/2}\partial
 h\|_{L^{\infty}(S_T)}
 \|r^{-1/4}\langle r\rangle^{-1/4}\nabla
 \phi\|_{L^2(S_T)}
\bigr)
\end{array}
\nonumber \\
 & & + C\log (2+T) \|{\bar \phi}_k\|^2_{Z_1(T-2/k)}
\nonumber \\
 & & \quad \times \bigl(\|h\|_{L^{\infty}(S_T)}
  + \|r^{1/2}\langle r \rangle^{1/2}\partial  h\|_{L^{\infty}(S_T)}
 \bigr)
 \nonumber
\end{eqnarray}
for $0 < t < T-2/k$.

To estimate $\|{\bar \phi}_k\|_{Z_1(T-2/k)}$ appearing on the
right-hand side of (\ref{egx3}),
we note that, by the Schwarz inequality and Lemma \ref{le1},
\begin{eqnarray*}
 \lefteqn{
\|r^{-1/4}\langle r \rangle^{-1/4}\partial
{\bar \phi}_k\|^2_{L^2(S_{T_k})}}\\
&=& \int_0^{T-2/k}\!\!\!\! \int_{{\mathbb R}^3} |x|^{-1/2}
\langle x\rangle^{-1/2} \\
& & \qquad \times\left|\int_0^T\!\! \int_{{\mathbb R}^3}
{\bar \rho}_k(t-s,x-y)\partial \phi(s,y)\,dyds\right|^2dxdt
\nonumber \\
&\leq& \int_0^T\!\! \int_{{\mathbb R}^3} |x|^{-1/2}
\langle x\rangle^{-1/2}
\int_{{\mathbb R}^3}
\rho_k(x-y)|\partial \phi(s,y)|^2\,dydxds
\nonumber \\
&\leq& C\|r^{-1/4}\langle r \rangle^{-1/4}\partial
\phi\|^2_{L^2(S_T)}.
\end{eqnarray*}
Similarly, we also have
\[
 \|r^{-5/4}\langle r \rangle^{-1/4}{\bar \phi}_k\|_{L^2(S_{T_k})}
 \leq C \|r^{-5/4}\langle r \rangle^{-1/4} \phi\|_{L^2(S_T)}.
\]
Hence, (\ref{egx3}) leads to
\begin{eqnarray}
\label{eg23}
\lefteqn{
\|\partial {\bar \phi}_k(t,\cdot)\|^2_{L^2({\mathbb R}^3)}
+\|{\bar \phi}_k\|^2_{Y_1(T-2/k)} + \|{\bar \phi}_k\|^2_{Z_1(T-2/k)}
} \\
&\leq& C\bigl(\|\nabla {\bar \phi}_k(0,\cdot)\|^2_{L^2} +
 \|\partial_t {\bar \phi}_k(0,\cdot)\|^2_{L^2}\bigr) \nonumber \\
 & & + C \bigl(\log (2+T)\bigr)^{1/2}
  \|\phi\|_{Z_1(T)}\|r^{1/4}\langle r \rangle^{1/4}F\|_{L^2(S_T)} \nonumber \\
 & & + C\log (2+T) \|\phi\|^2_{Z_1(T)}
\begin{array}[t]{l}
 \bigl(\|h\|_{L^{\infty}(S_T)} \\
 \quad + \|r^{1/2}\langle r \rangle^{1/2}\partial  h\|_{L^{\infty}(S_T)}
 \bigr).
\end{array}
 \nonumber
\end{eqnarray}

Finally, the inequality above passes to the limit as follows.

\bigskip

$\bullet$ $\displaystyle
\lim_{k \to \infty} \|\partial_{\alpha} {\bar
\phi}_k(t,\cdot)\|_{L^2({\mathbb R}^3)} = \|\partial_{\alpha}
\phi(t,\cdot)\|_{L^2({\mathbb R}^3)}
$ for $0 \leq t < T$.

\medskip

We may assume $0 \leq t < T - 2/k$. Then we have
\begin{eqnarray*}
 \partial_{\alpha} {\bar \phi}_k (t,x)
&=& \int_0^T \int_{{\mathbb R}^3} {\bar \rho}_k(t-s,
x-y)\partial_{\alpha} \phi(s,y) dyds, 
\end{eqnarray*}
and hence
\begin{eqnarray*}
 \lefteqn{\|\partial_{\alpha} {\bar \phi}_k(t,\cdot) - \partial_{\alpha} \phi(t,\cdot)\|_{L^2}}
 \\
&\leq& \int_0^T \int_{{\mathbb R}^3} {\bar \rho}_k(t-s,y)
\|\partial_\alpha \phi(s,\cdot - y) - \partial_\alpha \phi (t,\cdot)\|_{L^2}dyds
 \\
&\leq& \int_0^T \rho^{(1)}_k(t-s+1/k)
\|\partial_\alpha \phi(s,\cdot) - \partial_\alpha \phi (t,\cdot)\|_{L^2}ds
 \\
& & + \int_{{\mathbb R}^3} \rho_k(y)
\|\partial_\alpha \phi(t,\cdot - y) - \partial_\alpha \phi (t,\cdot)\|_{L^2}dy.
\end{eqnarray*}
Thus we have $
\lim_{k \to \infty} \partial_{\alpha} {\bar
\phi}_k(t,\cdot) = \partial_{\alpha} \phi(t,\cdot)$ in $L^2({\mathbb R}^3)$.

\medskip

$\bullet$ $\displaystyle \|\phi\|^2_{Y_1(T)} \leq \liminf_{k \to \infty}
\|{\bar \phi}_k\|^2_{Y_1(T-2/k)}$,\quad
$\displaystyle \|\phi\|^2_{Z_1(T)} \leq \liminf_{k \to \infty}
\|{\bar \phi}_k\|^2_{Z_1(T-2/k)}$.

\medskip

By definition,
\begin{eqnarray*}
 \lefteqn{(1+T-2/k)^{1/2}\|{\bar \phi}\|^2_{Y_1(T-2/k)}}\\
 &=& \int_0^{T-2/k}\!\!\!\!\int_{{\mathbb R}^3} r^{-5/2}|{\bar \phi}_k|^2 dxdt
  + \int_0^{T-2/k}\!\!\!\!\int_{{\mathbb R}^3} r^{-1/2}|\partial {\bar \phi}_k|^2 dxdt \\
 &=& \int_0^{T-2/k}\!\!\!\!\int_{{\mathbb R}^3} r^{-5/2}|{\bar \rho}_k {\bar *}
  {\bar \phi}|^2 dxdt
  + \int_0^{T-2/k}\!\!\!\!\int_{{\mathbb R}^3} r^{-1/2}|{\bar \rho}_k
  {\bar *} \overline{\partial \phi}|^2 dxdt.
\end{eqnarray*}
Hence we have $\|\phi\|^2_{Y_1(T)} \leq \liminf_{k \to \infty}
\|{\bar \phi}_k\|^2_{Y_1(T-2/k)}$ by Fatou's lemma.
Similarly, $\|\phi\|^2_{Z_1(T)} \leq \liminf_{k \to \infty}
\|{\bar \phi}_k\|^2_{Z_1(T-2/k)}$.
\hfill $\Box$

\subsection{Continuous Dependence on the Data}


\begin{lemma}
 \label{lg3}
Let $\phi \in X_2(T) \cap Y_2(T)$ be a solution of
 {\rm (\ref{ea1})}, {\rm (\ref{ea2})}. If
$\|\nabla f\|_{H^1({\mathbb R}^3)} + \|g\|_{H^1({\mathbb R}^3)} =
 \varepsilon \leq \varepsilon_1$
and $T \leq \exp (A_1/\varepsilon)$,
then we have
\begin{equation}
 \label{eg27}
 \|\phi\|_{E_2(T)} + \|\phi\|_{Y_2(T)} + \|\phi\|_{Z_2(T)} \leq M_1\varepsilon.
\end{equation}
Here, $\varepsilon_1$, $A_1$ and $M_1$ are the same constants as those in
 Lemma {\rm \ref{ld1}}.
\end{lemma}

Proof. Let ${\tilde T} \leq T$
be the largest number such that (\ref{eg27}) holds.
It is obvious that ${\tilde T} > 0$. Since it is now evident
from Lemma \ref{lg2} that
Lemma \ref{ld6} holds for $\phi, {\tilde \phi} \in X_2(T)$,
we get
\begin{eqnarray}
 \label{eg28}
\lefteqn{\|\phi\|_{E_2({\tilde T})}
+ \|\phi\|_{Y_2({\tilde T})} + \|\phi\|_{Z_2({\tilde T})}} \\
&\leq& C_3(\varepsilon + \varepsilon^2) + 2C_3 M_1\varepsilon
\|\phi\|_{Z_2({\tilde T})} \log (2+{\tilde T})
\nonumber
\end{eqnarray}
by considering (\ref{edx1}) in $S_{\tilde T}$ and letting
${\tilde \phi} := \phi$ in (\ref {edx2}).
Assume ${\tilde T} < T$. Since $T \leq \exp (A_1/\varepsilon)$,
it follows from (\ref{eg28}) and (\ref{ed43}) that
\[
 \|\phi\|_{E_2({\tilde T})}
+ \|\phi\|_{Y_2({\tilde T})} + \|\phi\|_{Z_2({\tilde T})}
< \frac{1}{2}M_1\varepsilon
+ \frac{1}{2}\|\phi\|_{Z_2({\tilde T})}.
\]
This gives $\|\phi\|_{E_2({\tilde T})}
+ \|\phi\|_{Y_2({\tilde T})} + \|\phi\|_{Z_2({\tilde T})} < M_1\varepsilon$,
which contradicts the definition of ${\tilde T}$.
Thus ${\tilde T} = T$. \hfill $\Box$

\begin{lemma}
 \label{lg5}
Let $\phi_1, \phi_2 \in X_2(T) \cap Y_2(T)$ be solutions of
{\rm (\ref{ea1})}, {\rm (\ref{ea2})}.
If $\|\nabla f\|_{H^1} + \|g\|_{H^1} = \varepsilon \leq \varepsilon_2$
and $T \leq \exp(A_1/\varepsilon)$,
we have $\phi_1 = \phi_2$.
\end{lemma}

Proof. It is obvious from Lemma \ref{lg2} that
we can apply Lemma \ref{ld9} to $\phi_1, \phi_2$.
By Lemma \ref{ld9} and Lemma \ref{lg3}, we have
\begin{eqnarray}
\lefteqn{\|\phi_1 - \phi_2\|_{E_1(1)} +
\|\phi_1 - \phi_2\|_{Y_1(1)}}
\nonumber \\
&\leq& C_4 \cdot 2M_1\varepsilon \cdot \sqrt{2} \cdot
\bigl(\|\phi_1 - \phi_2\|_{E_1(1)} +
\|\phi_1 - \phi_2\|_{Y_1(1)}\bigr).
\nonumber
\end{eqnarray}
Hence, noting (\ref{ed54}),
we have $\phi_1 = \phi_2$ for $\varepsilon \leq \varepsilon_2$ and $0 \leq t \leq 1$.
Applying the same estimate to $\phi_1(t+1,x)$ and $\phi_2(t+1,x)$,
we see that $\phi_1 = \phi_2$ for $1 \leq t \leq 2$.
Repeating this procedure, we conclude that $\phi_1 = \phi_2$
as long as $t \leq \exp(A_1/\varepsilon)$.
\hfill $\Box$

\bigskip

Proof of Theorem \ref{ta2}.
It is evident from Theorem \ref{ta1} and Lemma \ref{lg5} that the map $\Phi$
is well-defined.
Let $(f,g),\ ({\tilde f}, {\tilde g}) \in D(\varepsilon_0)$.
By Lemma \ref{ld9} and Lemma \ref{lg3}, we have
\begin{eqnarray}
\lefteqn{\|\Phi(f,g) - \Phi({\tilde f},{\tilde g})\|_{E_1(t)} +
\|\Phi(f,g) - \Phi({\tilde f},{\tilde g})\|_{Y_1(t)}}
\nonumber \\
&\leq& C_4\|(\nabla f, g) - (\nabla {\tilde f}, {\tilde g})\|_{L^2}
 \nonumber \\
& & + 2C_4 \cdot M_1\varepsilon_0 \cdot (1+t)^{1/2}
 \nonumber \\
& & \times \bigl(\|\Phi(f,g) - \Phi({\tilde f},{\tilde g})\|_{E_1(t)} +
\|\Phi(f,g) - \Phi({\tilde f},{\tilde g})\|_{Y_1(t)}\bigr).
\nonumber
\end{eqnarray}
Thus for $t \leq A_2\varepsilon_0^{-2}$,
\begin{eqnarray}
& &\|\Phi(f,g) - \Phi({\tilde f},{\tilde g})\|_{E_1(t)}
 + \|\Phi(f,g) - \Phi({\tilde f},{\tilde g})\|_{Y_1(t)}
\nonumber \\
& &\quad \leq 2C_4\left(
\|\nabla (f - {\tilde f})\|_{L^2} + \|g - {\tilde g}\|_{L^2}
\right). \nonumber
\end{eqnarray}
If we choose $N$ large enough so that
$N\cdot A_2\varepsilon_0^{-2} \geq T_{\varepsilon_0}$,
then we get
\begin{eqnarray}
& &\|\Phi(f,g) - \Phi({\tilde f},{\tilde g})\|_{E_1(T_{\varepsilon_0})}
 + \|\Phi(f,g) - \Phi({\tilde f},{\tilde g})\|_{Y_1(T_{\varepsilon_0})}
\nonumber \\
& &\quad \leq (2C_4)^N\left(
\|\nabla (f - {\tilde f})\|_{L^2} + \|g - {\tilde g}\|_{L^2}
\right). \nonumber
\end{eqnarray}
This completes the proof of Theorem \ref{ta2}.
\hfill $\Box$

\bigskip

{\bf Acknowledgements.}
The authors are very grateful to Jason Metcalfe for his significant suggestions.
They also thank the referees for providing constructive comments and criticisms,
which have improved the quality of this article.

\bibliographystyle{amsplain}


\end{document}